\documentclass[a4paper,12pt]{article}
\usepackage{geometry,times}
\usepackage{amsfonts}
\usepackage{amssymb}
\usepackage{amsmath}
\usepackage{amsthm}
\newcommand{\rose}{R\accent 23 u\v zi\v cka}
\newcommand{\bb}{\mathbb}

\newcommand{\intq}{\int\limits}
\theoremstyle{plain}
\newtheorem{satz}{Proposition}

\newtheorem{lemmab}[satz]{Lemma}
\theoremstyle{thm}
\newtheorem{theorem}{Theorem}
\theoremstyle{note}
\newtheorem{note}{Note}
\newtheorem{remark}[note]{Remark}
\theoremstyle{definition}
\newtheorem{defi}{Definition}
\theoremstyle{problem}
\newtheorem{problem}{Problem}
\theoremstyle{notadef}
\newtheorem{notationb}{Notation}
\usepackage[]{graphicx}
\usepackage{color}

\definecolor{ered}{rgb}{0.6,0,0}

\definecolor{dark-BLGR}{rgb}{0.28,0.28,0.28}

\definecolor{dark-red}{rgb}{0.90,0.14,0.14}

\definecolor{dark-viol}{rgb}{0.73,0.18,0.69}
\definecolor{Cyan}{rgb}{0.31,0.67,0.82}

\definecolor{light-yellow}{rgb}{1,1,0.8}

\definecolor{Blue}{rgb}{0.0353,0.0275,0.4}

\definecolor{GRAY}{rgb}{0.26,0.26,0.26}

\definecolor{GREYY}{rgb}{0.45,0.45,0.45}

\definecolor{vivid-viol}{rgb}{0.3255,0.0353,0.55}

\definecolor{dark-blue}{rgb}{0.05,0.05,0.65}

\definecolor{dark-green}{rgb}{0.03,0.77,0.29}

\definecolor{dark-Green}{rgb}{0.03,0.57,0.09}

\definecolor{strong-viol}{rgb}{0.2353,0.094,0.349}

\definecolor{BLUE}{rgb}{0.41,0.44,0.93}

\definecolor{RED}{rgb}{0.90,0.18,0.32}

\definecolor{Yel}{rgb}{0.89,0.95,0.19}

\definecolor{White}{rgb}{1,1,1}
\definecolor{black}{rgb}{0,0,0}

\definecolor{Orchid}{rgb}{0.6,0.1607,0.2}

\definecolor{Orange}{rgb}{0.99,0.49,0}

\definecolor{Thistle}{rgb}{0.1151,0.1249,0.9873}

\definecolor{brown}{rgb}{0.298,0.153,0.07843}

 \usepackage{titlesec}
\titleformat{\section}
{\normalfont
  \bfseries
}
{\thesection}{1ex}{}
\titleformat{\subsection}
{\normalfont
  \bfseries
}
    {\thesubsection}{1ex}{}
\begin{document}
\hspace*{1cm}\\[2cm]
\textbf{The Stokes Eigenvalue Problem on balls and annuli in three
  dimensions: Solutions  with Poloidal and Toroidal Fields}
\\[.5cm]
Bernd Rummler\\
\hspace*{2mm}{\small \em Otto-von-Guericke-Universit\"at Magdeburg,
  Inst. f\"ur Analysis und Numerik, PF
  4120, 39016 Magdeburg}
\\[3mm]
Gudrun Th\"ater\footnote{Corresponding author\quad
  E-mail:~\textsf{gudrun.thaeter@kit.edu}}\\
\hspace*{2mm}{\small \em KIT, Inst. f\"ur Angewandte \& Numerische Mathematik,
  76128 Karlsruhe} 
\\[5mm]
\hspace*{.1cm}\hfill\parbox{16cm}{
{\small\textbf{Abstract:}
We consider the Stokes eigenvalue problem in open balls and open
annuli in ${\bb R}^{3}$ with homogeneous Dirichlet boundary
conditions. Using the frame of toroidal and poloidal fields we
construct the othogonal decomposition of the Stokes eigenvalue problem
in problems for toroidal and poloidal eigenfunctions.  This provides
the proof of the completeness of a system of explicitly calculated Stokes eigenfunctions 
given by one of the authors in 1999, \cite{RumHab} }.}
\section{Introduction\label{sec_int}}
Flow problems in open balls and open annuli play an important role in
applications, \cite{RRTZAMM}. Due to their very
symmetric geometry, it is possible
to gain much more precise information than in general 
domains. This helps to solve those flow problems analytically as well
es numerically.

In a series on {\em Eigenfunctions
of the Stokes Operator in Special Domains} \cite{LeeRum, RumSto1,
RumSto2}, we have constructed (i.e. calculated) 
complete systems of (complex-valued) eigenfunctions for those
domains. Proofs of completeness of these systems and the
real-valued eigenfunctions were
included in \cite{RumHab}.

The paper at hand returns to these questions with a different goal:
We establish a general method to decompose solenoidal vector fields on a ball
and a spherical annulus (a spherical layer). They can also be used to
prove completeness for 
the system of calculated real-valued eigenfunction  
written in Cartesian coordinates, \cite{RumHab}.

It is worth to note, that in \cite{RumHab} we have implicitly
decomposed solenoidal fields  
into toroidal and poloidal vector fields. However -- in contrast to
here --
there we have only used 
the shape of toroidal and poloidal vector fields (written in Cartesian 
coordinates as functions of spherical coordinates) in order to exploit the
scalar action of the vector Laplacian on every  
Cartesian coordinate.

The decomposition of solenoidal fields
into toroidal and poloidal vector parts was introduced by Gustav Mie in 1908
\cite{Mie} in  the mo\-deling of optical properties of colloidal
metallic solutions and is sometimes called the Mie re\-presentation
of solenoidal fields. 
Since 1957,  the Mie re\-presentation has been mainly applied
in geomagnetism. In particular, there the decomposition is used 
to analyze and to interpret  measurement data of Earth observation satellites.
Moreover, since nearly thirty years it has been applied 
in fluid dynamics as well for, e.g.,  the study of
the Boussinesq equations.

It is worth mentioning that one can understand the decomposition of solenoidal fields 
into toroidal and poloidal parts as a refinement of the 
standard Helmholtz decomposition of three-dimensional (3D) vector fields
in such a way, that one divides the set of solenoidal fields 
into two construction types. But 
the method  is restricted to very special domains, e.g.
such domains in $ {\bb R}^{3}$, where
the prescribed Dirichlet boundary conditions on the rigid
part of the boundary can be attached to a unique normal direction
(a unique normal unit vector). The decomposition is possible for all
eigenvalue problems  
of the Stokes operator in 3D domains, which allow the explicite determination 
of the eigenfunctions by the method of separation of variables.
From this point of view one can reach the results of \cite{LeeRum,
  RumSto1, RumSto2} also
by the use of decompositions into toroidal and poloidal vector fields.

We examine the Stokes eigenvalue problem on two
types of domains: the open 
spherical annulus ${\Omega}_{\sigma}$ (the gap between two 
concentrical spherical surfaces, where the radius $r$ runs from a 
fixed $0\,<\, \sigma\,<\, 1$ to $1$)
as well as  the open  
unit ball ${\Omega}_{o}$ (since ${\sigma}=0$)
both in $ {\bb R}^{3}$. 
Before we can turn to our primary objective which is to
prove the promised decomposition theorem 
we start 
with the formulation of the Stokes eigenvalue problem (for the
detailed notation
see Subsection 2.1):\\[-.5cm]
\begin{problem} \label{Prob1} We seek solutions ${\underline{u}}$, $\lambda$ and $p$
(for ${\sigma}:0\leq{\sigma}<1$)  fulfilling:\\[.2cm]
$- \triangle {{\;\!}}{\underline{u}}\,
 +\,{\underline{\nabla}} \hspace*{.1cm} p
 = \,\lambda {\underline{u}}\,\quad$ and
$\,\quad div \hspace*{.1cm}{\underline{u}}\,=\,
{\underline{\nabla}}^{T}\cdot {\underline{u}}\,=\,0\quad$
{in} ${\Omega}_{\sigma}$ , with 
${\underline{u}}\,=\,{\underline{0}}$ on $\partial{\Omega}_{\sigma}$\,.
\end{problem}
{~}\\[-.2cm]
The operator formulation in appropriate function spaces of
solenoidal vector fields  \cite{CoFoi, RumHab, Temam}
reads:{~}\\[-.5cm]
\begin{problem} \label{Prob2} Find ${\underline{u}}\in
  D({\bf A}_{\sigma})$ and $\lambda$ fulfilling: ${\bf
    A}_{\sigma}{\underline{u}}\,=\,\lambda {\underline{u}}$ , where
  ${\bf A}_{\sigma}$ denotes the Stokes operator.
\end{problem}
{~}\\[-.5cm]
The operator
${\bf A}_{\sigma}$ (for fixed ${\sigma}:0\leq{\sigma}<1$) is
positive self-adjoint with pure real point spectrum, 
i.e., all eigen\-values $\lambda_{j}\,>\,0$ are of 
finite multiplicity. The associated eigenfunctions
$\{{\underline{w}}_{j}({\underline{x}})\}_{j=1}^{\infty}$ of 
${\bf A}_{\sigma}$ (counted in the multiplicity of $\lambda_{j}$) are
re\-gular, \cite{ADN, Catbri, Temam, Triebel}.
Problems 1 and 2  are well defined in the classical sense, i.e. 
with smooth ${\underline{\cal C}}^{m}$-eigenfunctions
$\{{\underline{w}}_{j}({\underline{x}})\}_{j=1}^{\infty}$ ($m\geq 2$).
The  justification is provided by the choice of homogeneous Dirichlet
boundary conditions on the boundary
$\partial{\Omega}_{\sigma}$ of ${\Omega}_{\sigma}$ in the
definition of the Stokes operator (details in Subsection
\ref{sec_spac_not}).

We construct toroidal and 
poloidal vector fields as special fields
defined by toroidal and poloidal scalar functions (potentials) respectively, 
to emphasize the structure of the eigenfunctions
$\{{\underline{w}}_{j}({\underline{x}})\}_{j=1}^{\infty}$.  
Espe\-cially we use the nice properties of the Lapace-Beltrami
operator ${\bf B}$ on spherical surfaces 
and the surface spherical harmonics as eigenfunctions of ${\bf B}$
for the proof of the decomposition theorem. 
The uniqueness of the toroidal and poloidal parts of every solenoidal vector field
is ensured by a vanishing mean-value conditions on spherical surfaces.
\\[3mm]
This paper is organized as follows:
All notation such as symbols, operators and function spaces and 
in particular the definition of the
Stokes operator ${\bf{A}}_{\sigma}$  
and the Lapace-Beltrami operator ${\bf B}$ is introduced in Section 2.
In Section 3 we provide an essential property for the characterization of solenoidal 
fields  on $\Omega_{\sigma}$ as Lemma
\ref{LE1} and formulate our decomposition Theorem \ref{Dec_Thm} in
Subsection \ref{sec_thm_dec}. 
There the real-valued eigenfunctions of ${\bf B}$:
$\{Z_{l}^{k}\}_{k=-l}^{l}$,  $\forall \, l=0,1,2,\dots\,$  are used 
as a technically efficient tool in the proof of Theorem
\ref{Dec_Thm} (cf. {{Theorem B}} in the Appendix).

The Stokes eigenvalue problem is investigated in Section 4. The
decomposition into problems for toroidal and poloidal eigenfunctions
especially using eigenvalue problems for the  
constitutive scalar functions are provided. There we derive the Stokes
ope\-rator in the 
unknowns $\lambda$ and $\underline{u}$ (the eigenvalue and the
toroidal or poloidal vectorial 
eigenfunction) from the original eigenvalue problem in the unknowns
$\lambda$ and ${\psi}$ resp. ${\chi}$ 
(the eigenvalue and a scalar eigenfunction).

Subsequently the eigenvalue problem of the Stokes operator
${\bf{A}}_{o}$ is considered on the open  
unit ball ${\Omega}_{o}$ in ${\bb R}^{3}$ with homogeneous Dirichlet boundary
conditions. We note that 
the decomposition of the Stokes eigenvalue problem in problems for
toroidal and poloidal eigenfunctions 
leads also to the proof of the completeness of a system of explicitly
calculated Stokes eigenfunctions from \cite{RumHab}.

Additionally in Subsection \ref{sec_sto_ewAnn} we show the general
shape of the eigenfunctions and the governing trans\-cendental  
equations for the eigenvalues for the Stokes ope\-rator ${\bf{A}}_{\sigma}$ on  
open spherical annuli  ${\Omega}_{\sigma}$.
The proof of completeness of these systems of explicitly calculated eigenfunctions
is quite similar to the proof for the unit ball.

\section{Notation and Basics\label{sec_not}}
\subsection{Notation for vector fields\label{sec_nov_not}}
Let ${\bb N}$, ${\bb N}_{o}$, ${\bb R}$, and ${\bb C}$ be the sets of 
positive integers, non-negative integers, real numbers and complex
numbers, respectively.  
For $n\,\in\,{\bb N}$ we denote by ${\bb R}^{n}$ the (real)
$n$-dimensional Eu\-cli\-dean space. 
We write underlined small letters ${\underline{x}}$,
${\underline{y}}$, ${\underline{u}}$, ${\underline{v}}$, $\dots$ 
to represent elements of ${\bb R}^{n}$. We understand 
these elements in the classical sense of linear algebra as columns,
i.e. 
${\underline{x}}\,=\,(x_{1},x_{2},\dots ,x_{n})^{T}$
and the symbol ${~}^{T}$ is transposition.
\\ [.3cm]
The space ${\bb R}^{n}$ is a  Hilbert space when 
equipped with the Eucli\-dean inner product:
$({\underline{x}},{\underline{y}})_{\cal{E}}\,=\, {\underline{x}}^{T}{\underline{y}}\, 
$ $\forall\, {\underline{x}}\,,\,{\underline{y}}\,\in\,{\bb R}^{n}$
and the Euclidean norm
$\|{\underline{x}}\|_{\cal{E}}\,=\,
\sqrt{{\underline{x}}^{T}{\underline{x}}}\, $,
respectively. 

We will assume, that ${\Omega}$ is any 
open bounded set of class ${\cal C}^{2}$ in ${\bb R}^{n}$ with 
$n\,\in\,{\bb N}$. We write ${\overline{\Omega}}$ for the
closure of the set ${\Omega}$ and $\partial{\Omega}\,=\,{\overline{\Omega}}\setminus{\Omega}$
for the boundary of ${\Omega}$.
Vector fields ${\underline{u}}:$ $D({\underline{u}}) \subset {\bb R}^{n}\,\rightarrow \,{\bb R}^{n}$, 
 with $D({\underline{u}})\,\in\, \{{\Omega},{\overline{\Omega}},\partial{\Omega}\}$,
are regarded  as  functions of the vector of positions
${\underline{x}}$. 
It is usual to call the vectorial differential expression
\begin{equation} \label{nab..def:}
{\underline{\nabla}}\,:=\,(
\frac{\displaystyle{\partial}}{\displaystyle{\partial x_{1}}},
\frac{\displaystyle{\partial}}{\displaystyle{\partial x_{2}}},\dots,
\frac{\displaystyle{\partial}}{\displaystyle{\partial x_{n}}})^{T}
\end{equation}
the Nabla operator (here ${\underline{\nabla}}$ is written in canonical 
(Cartesian) coordinates).

Our considerations will be focused on problems, where $n=3$ (or $n=2$ for 
tangential components on surfaces).
In particular, we consider the open ball ${\Omega}:={\Omega}_{o}$ and
the open annulus ${\Omega}:={\Omega}_{\sigma}$, $0<{\sigma}<1$.

\begin{notationb}\label{N1}
We write
${\Omega}_{o}\,:=\{{\underline{x}} \in {\bb R}^{3}:\,
\|{\underline{x}}\|_{\cal{E}}\,
<1\}$  for the open unit ball
and
$\omega:=\partial{\Omega}_{o}\,
:=\{{\underline{x}} \in {\bb R}^{3}:\,\|{\underline{x}}\|_{\cal{E}}=1\}\,$ for 
its  boundary.
The open spherical layer (spherical annulus) 
will be denoted by
${\Omega}_{\sigma}\,:=\{{\underline{x}} \in {\bb
  R}^{3}:\,{\sigma}<\|{\underline{x}}\|_{\cal{E}}<1\}$ and its
boundary is
$\partial{\Omega}_{\sigma}\,=\,
\omega_{\sigma} \cup \omega$.
The spherical surface $\omega_{r}$ is defined by
$\omega_{r}\,:=\,\{{\underline{x}}\in {\bb
  R}^{3}:\,\|{\underline{x}}\|_{\cal{E}}=r\} 
\,=\,r\omega$
for all $r\,\in\,(0,1]$.
Throughout this paper we will use $\omega\,=\,\omega_{1}$ 
as the notation for the unit spherical surface in ${\bb R}^{3}$.
\end{notationb}
\begin{notationb}\label{N2} We write twofold  
underlined capital letters ${\underline{\underline{A}}}\,$ , ${\underline{\underline{B}}}\,$ , $\dots$
to represent
(3,3)-matrices.
For 
${\underline{\underline{A}}}$ and
${\underline{\underline{B}}}$ we use the so-called Frobenius inner product: 
\begin{equation} \label{17..def:}
\hspace*{3cm}{\underline{\underline{A}}}:{\underline{\underline{B}}}:=
tr({\underline{\underline{A}}}\cdot{\underline{\underline{B}}}^{T})\,=\,
\sum_{j=1}^{3}\sum_{k=1}^{3}a_{j,k}b_{j,k},
\end{equation}
where $tr(.)$ denotes the trace of the (product-)matrix.
\end{notationb}
\noindent The vector product of  ${\underline{u}}\,\,\mbox{and}\,{\underline{v}}\,$
as elements of ${\bb R}^{3}$ is defined  by
${\underline{u}}\,\times\,  {\underline{v}}$. It 
is ortho\-gonal to 
${\underline{u}}({\underline{x}})\,\mbox{and}\,{\underline{v}}({\underline{x}})$
in the sense of ${\bb R}^{3}$ at any the point
${\underline{x}}\in\,D({\underline{u}})\cup D({\underline{v}})$.
\\[.3cm]
To distinguish the representation of vector fields (or vector-valued 
differential expressions) in Cartesian and 
especially spherical polar coordinates we denote 
by ${\underline{u}}_{\mathfrak{c}}$ the vector field
${\underline{u}}$ written in a Cartesian coordinate system 
and by ${\underline{u}}_{\mathfrak{s}}$ the field ${\underline{u}}$ written in
a polar coordinate system (or e.g. by ${\underline{\nabla}}_{\mathfrak{s}}$ the Nabla operator in spherical polar coordinates
).

\begin{notationb}\label{N3} Let the unit vectors in the Cartesian coordinate system
be given by ${\underline{\mathfrak{e}}}_{j}\,:=\,
(\delta_{j,1},\delta_{j,2},\delta_{j,3})^{T}$  $\forall\,j=1,2,3$, with
Kronecker's delta $\delta_{j,k}$.
If we denote by $\{{\underline{\mathfrak{e}}}_{r},
{\underline{\mathfrak{e}}}_{\vartheta},{\underline{\mathfrak{e}}}_{\varphi}\}$
the system of unit vectors in spherical polar coordinates then ${\underline{u}}$
is representable in both systems as\\[.1cm]
\hspace*{2.3cm}${\underline{u}}\,=\,\sum_{j=1}^{3} u_{j}{\underline{\mathfrak{e}}}_{j}\,=\, \sum_{j=1}^{3} u_{j,{\mathfrak{c}}}{\underline{\mathfrak{e}}}_{j}\,=\, 
 u_{r}{\underline{\mathfrak{e}}}_{r}+
 u_{\vartheta}{\underline{\mathfrak{e}}}_{\vartheta}+
 u_{\varphi}{\underline{\mathfrak{e}}}_{\varphi}\,=\, 
 u_{r,{\mathfrak{s}}}{\underline{\mathfrak{e}}}_{r}+
 u_{\vartheta,{\mathfrak{s}}}{\underline{\mathfrak{e}}}_{\vartheta}+
 u_{\varphi,{\mathfrak{s}}}{\underline{\mathfrak{e}}}_{\varphi}$.
\end{notationb}
\begin{remark}\label{R1} The transformation from one coordinate system to the other 
can be written as ${\underline{u}}_{\mathfrak{c}}\,=\,
{\underline{\underline{T}}}_{{\mathfrak{c}},{\mathfrak{s}}}{\underline{u}}_{\mathfrak{s}}$
or ${\underline{u}}_{\mathfrak{s}}\,=\,
{\underline{\underline{T}}}^{-1}_{{\mathfrak{c}},{\mathfrak{s}}}
{\underline{u}}_{\mathfrak{c}}\,=\,
{\underline{\underline{T}}}_{{\mathfrak{s}},{\mathfrak{c}}}{\underline{u}}_{\mathfrak{c}}$, respectively,
where we have used the concept of columns of coordinates
and the transformation matrices
\begin{equation*}
		{\underline{\underline{T}}}_{{\mathfrak{c}},{\mathfrak{s}}}\,
		:=\,\left[
		\begin{array}{lll}
\sin{\vartheta}\cos{\varphi}& \cos{\vartheta}\cos{\varphi}& -\sin{\varphi} \\
\sin{\vartheta}\sin{\varphi}& \cos{\vartheta}\sin{\varphi}& {~}\cos{\varphi}\\	
\cos{\vartheta} & -\sin{\vartheta}& {~} 0
		\end{array}
		\right] \quad\,\mbox{and} \quad
{\underline{\underline{T}}}_{{\mathfrak{s}},{\mathfrak{c}}}\,:=\,
{\underline{\underline{T}}}_{{\mathfrak{c}},{\mathfrak{s}}}^{-1}\,=\,
{\underline{\underline{T}}}_{{\mathfrak{c}},{\mathfrak{s}}}^{T}\,.
\end{equation*}
\end{remark}
\begin{remark}\label{R2} 
We use the Nabla operator ${\underline{\nabla}}$ to represent
standard operations of  
vector analysis in ${\bb R}^{3}$ as
$$ div{\;\!}{\underline{u}}({\underline{x}})\,:=\,{\underline{\nabla}}^{T}{\underline{u}}({\underline{x}}), \quad
curl{\;\!}{\underline{u}}({\underline{x}})\,:=\,{\underline{\nabla}}{\;\!}\times{\;\!}
{\underline{u}}({\underline{x}}), \quad \mbox{and}\quad
grad{\;\!}{\;\!}{\eta}({\underline{x}})\,:=\,{\underline{\nabla}}{\;\!} {\eta}({\underline{x}})$$
for vector fields ${\underline{u}}$ and scalar functions ${\eta}$.
\end{remark}
\noindent To some extend we will use these 
as projections to define (surface-)gradient, curl, 
and divergence in Section 3 as well.

\subsection{Function spaces and auxiliarities\label{sec_spac_not}}
We regard 
the usual multi-indeces ${\underline{\kappa}}$ also as columns.
Following Schwartz' notation we use
${\underline{\kappa}}\,\in\,{\bb N}_{o}^{n}$ as a multi-index, with
the degree $|{\underline{\kappa}}|\,:=\,\sum_{j=1}^{n}{\kappa}_{j}$
and write partial derivatives as 
\begin{equation*}
D^{\underline{\kappa}}\,:=\,\frac{\displaystyle{{\partial}^{|{\underline{\kappa}}|}}}
{\displaystyle{{\partial}x_{1}^{{\kappa}_{1}}\dots{\partial}x_{n}^{{\kappa}_{n}}}}
\,\,,\,\forall \, {\underline{\kappa}} \in\,{\bb N}_{o}^{n}\,\quad\mbox{with}\,\quad
|{\underline{\kappa}}|\geq 1\,.
\end{equation*}
For ${\underline{\kappa}}\,=\,{\underline{0}}$ and any 
function
$f:D(f)\,\subset\,{\bb R}^{n}\rightarrow {\bb R} \,(\mbox{or}\,{\bb C})$
we set $D^{\underline{0}}f({\underline{x}})\,:=\,f({\underline{x}})
\,\,\forall \, {\underline{x}}\,\in\,D(f)$.
\begin{notationb}\label{N4}  We denote by ${\bb C}^{m}({\overline{\Omega}})$,  $m\,\in\,{\bb N}$,
the Banach space of real-valued 
(or complex-valued) functions $f$ that are continuous on ${\overline{\Omega}}$
and all of whose derivatives $D^{\underline{\kappa}}f$ up to and 
including the order $m=|{\underline{\kappa}}|$ are continuous
on ${\Omega}$ and can be extended by continuity to ${\overline{\Omega}}$.
We write ${\bb C}({\overline{\Omega}})$ for $m\,=\,0$.
\end{notationb}
\begin{notationb}\label{N5}  Let ${\tilde{\Omega}}$ be any simply connected non-empty set in ${\bb R}^{n}$.
With 
${\cal C}^{\infty}({\tilde{\Omega}})$ we denote 
the linear space of all real-(or complex-)valued functions 
$f$ which in all points ${\underline{x}}$  
of  ${\tilde{\Omega}}$ in ${\bb R}^{n}$ have continuous derivatives of
arbitrary order $|{\underline{\kappa}}|\,\in\,{\bb N}_{o}$ 
and  ${\cal C}^{\infty}_{o}({\tilde{\Omega}})$ is the subspace of
${\cal C}^{\infty}({\tilde{\Omega}})$ 
of functions
with 
compact supports.
\end{notationb}
\begin{notationb}\label{N6}  The Lebesgue spaces ${\bb
    L}_{2}({\Omega})$  and 
${\bb L}_{2}(\partial{\Omega})$ 
are spaces of \it{(equivalence classes)} of absolutely quadratic
integrable functions w.r.t. to the Lebesgue measure.
The Sobolev spaces ${\bb W}_{2}^{m}({\Omega})$ and
${\bb W}_{2}^{m}({\partial{\Omega}})$ are subspaces of 
${\bb L}_{2}({\Omega})$ (or of ${\bb L}_{2}(\partial{\Omega})$) with
derivatives up to and 
including the order $m=|{\underline{\kappa}}|$ in 
${\bb L}_{2}({\Omega})$ (or in ${\bb L}_{2}(\partial{\Omega})$, respectively).
The Hilbert spaces ${\bb L}_{2}(.)\,=\,{\bb W}_{2}^{0}(.)$ and ${\bb W}_{2}^{m}(.)$ are 
equipped with 
the norm 
$$
\|f\|_{{\bb W}_{2}^{m}(.)}\,:=\,
(\sum_{m \leq
  |{\underline{\kappa}}|}(\|D^{\underline{\kappa}}f\|_{{\bb
    L}_{2}(.)})^{2})^{(1/2)}\,,\quad m\,\in\,{\bb N}_{o}\,.
$$
We explain the spaces ${\bb W}_{2,o}^{m}({\Omega})$
of functions with 
vanishing boundary values in the generalized sense  
for $m\,\in\,{\bb N}$, as the closure of
${\cal C}^{\infty}_{o}({\Omega})$ in ${\bb W}_{2}^{m}({\Omega})$.
Finally we write the Sobolev-Slobodeckij spaces of order
$m+1/2$ on $\partial{\Omega}$ as 
${\bb W}_{2}^{m+1/2}({\partial{\Omega}})$ for $m\,\in\,{\bb N}_{o}$. The dual space of 
${\bb W}_{2}^{1/2}
({\partial{\Omega}})$ is denoted by ${\bb
  W}_{2}^{-1/2}({\partial{\Omega}})$.
\end{notationb}
\begin{notationb}\label{N7}
For any linear space ${\cal X}$ or any Banach space ${\bb X}$ we will denote
by  ${\underline{\cal X}}$ resp. ${\underline{\bb X}}$ the spaces
of vector-valued functions with all components in  ${\cal X}$ resp. 
${\bb X}$.\end{notationb}
\begin{notationb}\label{N8}  The abbreviation $\langle {\overline{f}}
  \rangle({\underline{x}})$  is the radial mean value 
of functions $f\,\in\,{\bb L}_{1}(\tilde{\Omega})$ 
 \begin{equation}\label{DefRad}
\langle {\overline{f}}\rangle({\underline{x}})\,=\,
\langle {\overline{f}}\rangle
(r)\,
\:=\,\frac{1}{4\pi}{\intq_{{\underline{y}}\in {\omega}}{\;\!}f(r \cdot {\underline{y}})\,d {\omega}}\quad \mbox{with}\quad  
\,r=\|{\underline{x}}\|_{\mathcal{E}}\,\,,
\end{equation}
where $ \tilde{\Omega}$ is one of the rotational symmetrical
domains ${\Omega}_{o},\,\omega,\,{\Omega}_{\sigma},\,\omega_{\sigma}$
in ${\bb R}^{3}$. We note that $\langle {\overline{f}} \rangle(r)$ is well defined for almost all $r\in [{\sigma},1]$
by Fubini's theorem.
\end{notationb}
\begin{notationb}\label{N9}  
We will say that the vector fields ${\underline{u}}$ are {{\em solenoidal}} on  
${\Omega}_{o}$ or ${\Omega}_{\sigma}$, if  
$div{\;\!}{\underline{u}}\,=\,0$
there (at least in the generalized sense) and
for almost all $r>0: \,{\omega}_{r}\,\subset \,{\overline{\Omega}}_{o}
$ or $\,{\omega}_{r}\,\subset
\,{\overline{\Omega}}_{\sigma}$, respectively, it holds
 \begin{equation}\label{Defsolfield}
0\,=\,\intq_{{\omega}_{r}}{\;\!}{\underline{u}}^{T}{\underline{x}}\,d {\omega}_{r}  
\end{equation},
where we understand the integral at least in the sense of generalized
functions on  ${\omega}_{r}$ too.
\end{notationb}
\begin{notationb}\label{N10}  Let us define by 
${\underline{\cal V}}_{\sigma}\,:=\,
\{{\underline{v}}\,\in\,
{\underline{{\cal C}}}^{\infty}_{o}({\Omega_{\sigma}}):div{\;\!}{\underline{v}}=0\}$
for $ {\sigma}:0\leq{\sigma}<1$. The closures of ${\underline{\cal V}}_{\sigma}$
in the sense of ${\underline{\bb L}}_{2}({\Omega}_{\sigma})$ will be
called ${\underline{\bb S}}({\Omega}_{\sigma})$  
and the closures of ${\underline{\cal V}}_{\sigma}$ in
${\underline{\bb W}}_{2,o}^{1}({\Omega}_{\sigma})$  
will be called ${\underline{\bb S}}^{1}({\Omega}_{\sigma})$, 
respectively. We abbreviate the space ${\underline{\bb
    W}}_{2}^{2}({\Omega}_{\sigma})\cap{\underline{\bb
    S}}^{1}({\Omega}_{\sigma})$ 
with ${\underline{\bb S}}^{2}({\Omega}_{\sigma}):={\underline{\bb
    W}}_{2}^{2}({\Omega}_{\sigma})\cap{\underline{\bb
    S}}^{1}({\Omega}_{\sigma})$.
\end{notationb}
\begin{remark}\label{R3}
  The only purpose of condition (\ref{Defsolfield}) in  Notation \ref{N9} is it 
to exclude scalar multiples of the harmonic field
$grad{\;\!}({\|{\underline{x}}\|_{\cal{E}}})^{-1}$ at the outset. 
The advantage of the spaces introduced in Notation \ref{N9} is, that
their elements are solenoidal fields and fulfill the boundary conditions
in a generalized sense. Especially the justification of the existence
of the trace of the normal component for
${\underline{u}}\,\in\,{\underline{\bb S}}({\Omega}_{\sigma})$ in 
the ${\bb W}_{2}^{-1/2}({\partial{\Omega}})$ 
sense is given by the condition
$div{\;\!}{\underline{u}}\,\in\,{\underline{\bb
    L}}_{2}({\Omega}_{\sigma})$. 
One can understand solenoidal fields on
${\Omega}_{\sigma}\,,\,\sigma\,>\,0$, also in the sense of 
elements of the range: $curl{\;\;\!}{\underline{\bb
    W}}_{2}^{1}({\Omega}_{\sigma})$ (cf. \cite{Temam}
p.467).
\end{remark} 
\begin{defi} \label{D2}
  We denote by 
$\triangle \,=\, div{\;\!}grad{\;\!}:{\underline{\bb W}}_{2}^{2}({\Omega}_{\sigma})
\rightarrow{\underline{\bb L}}_{2}({\Omega}_{\sigma})$
the Laplace operator in the sense of Friedrichs extension.
The Leray-Helmholtz projector $\Upsilon$ is the projector
on solenoidal fields explained by: $\Upsilon : {\underline{\bb L}}_{2}({\Omega}_{\sigma})
\rightarrow{\underline{\bb S}}({\Omega}_{\sigma})$.
The Stokes operator is stated as the product of $\Upsilon$
and $-\triangle$ and is defined by ${\bf A}_{\sigma}\,:=\,- \Upsilon \triangle$, where
$- \Upsilon \triangle :{\underline{\bb S}}^{2}({\Omega}_{\sigma})
\rightarrow{\underline{\bb S}}({\Omega}_{\sigma})$.\end{defi}
\begin{remark}\label{StokesOp}
  The Stokes operator ${\bf A}_{\sigma}$ is a positive self-adjoint operator
with compact inverse. This follows in a standard way from the general theory 
of elliptical systems. It also provides 
the existence of a countable system of Stokes eigenfunctions
and the completeness of such a system in the subspace of the Hilbert space 
${\underline{\bb L}}_{2}(.)\,=\,({\bb L}_{2}(.))^{3}$ 
consisting of weak solenoidal vector fields which fulfill the generally 
required  boundary conditions in the ${\bb W}_{2}^{-1/2}$ trace sense.
The exact knowledge of the eigenfunctions of the Stokes operator is 
especially useful for the the numerical treatment of the incompressible
Navier-Stokes equations by Galerkin methods.
\end{remark}
\begin{defi} \label{D3} We denote by $(.,.)_{D}$ the Dirichlet (scalar) product
defined for vector fields\\[.1cm]
${\underline{u}}$ and ${\underline{v}}$
from ${\underline{\bb W}}_{2}^{1}({\Omega}_{\sigma})$ through:
\begin{eqnarray} \label{Dirichlet} 
({\underline{u}},{\underline{v}})_{D} \:\! := 
\intq_{{\Omega}_{\sigma}}{\:\!}(\sum_{j=1}^{3}\sum_{k=1}^{3}{\;\!}
{\frac{\displaystyle{\partial u_{j}}}{\displaystyle{\partial x_{k}}}}{\frac{\displaystyle{\partial v_{j}}}{\displaystyle{\partial x_{k}}}})
\,d {\Omega}_{\sigma}
 \:\! = \intq_{{\Omega}_{\sigma}}{\:\!}({\underline{\nabla}}\cdot{\underline{u}}^{T}({\underline{x}}))^{T}:
({\underline{\nabla}}\cdot{\underline{v}}^{T}({\underline{x}}))^{T}\,d {\Omega}_{\sigma}
\;\, \forall \:\! {\underline{u}} ,\:\! {\underline{v}}\:\!  \in \:\! {\underline{\bb W}}_{2}^{1}({\Omega}_{\sigma})\,,
\end{eqnarray}
The matrices are here the Jacobian matrices
${\underline{\underline{A}}}:=({\underline{\nabla}}\cdot{\underline{u}}^{T})^{T}$ and  ${\underline{\underline{B}}}:=({\underline{\nabla}}\cdot{\underline{v}}^{T})^{T}$. \end{defi}
\begin{remark}\label{R4}  The Dirichlet  (scalar) product $(.,.)_{D}$ generates an equivalent
norm on ${\underline{\bb W}}_{2,o}^{1}({\Omega}_{\sigma})$ and a seminorm on 
${\underline{\bb W}}_{2}^{1}({\Omega}_{\sigma})$.\end{remark}
\begin{defi} \label{D4}
The Laplace-Beltrami operator is defined as
\begin{equation*}
{\bf B^{o}}\,{Y}\, :=\, -\,
\frac{\displaystyle{1}}
{\displaystyle{\sin{\vartheta}}}
\frac{\partial }{\partial \vartheta}(\sin{\vartheta}\frac{\partial Y}{\partial \vartheta})
\,-\,
\frac{\displaystyle{1}}{\displaystyle{\sin ^{2}{\vartheta}}}
\frac{\displaystyle{\partial ^{2} Y}}
{\displaystyle{\partial \phi^{2}}}\,,\, \hspace*{0.7cm}  \forall\,Y\,\in\,D({\bf B^{o}})\,=\,C^{\infty}(\omega)\,\subset\,{\bb L}_{2}(\omega)\,.
\end{equation*}
We denote the Friedrichs' extension of ${\bf B^{o}}$ by 
${\bf B}:={\overline{\bf B^{o}}}$. The operator ${\bf B}$ is 
called Laplace-Beltrami operator as well. 
It is convenient to use the eigenvalues and eigenfunctions of the Laplace-Beltrami operator for the definion of 
function spaces on $\omega$ (cf. e.g.{\cite{CoFoi}}, p. 33). 
In this way, we introduce the Hilbert spaces 
$D({\bf B}^{\frac{m}{2}})$  for $m\in\,{\bb N}_{o} $ by\\[-.7cm]
\begin{eqnarray} \label{BeltrSpaces} 
D({\bf B}^{\frac{m}{2}})\,=\,\{Y\,\in\,{\bb L}_{2}(\omega):\,
\sum_{j=0}^{\infty}\,(1+\nu_{j}^{m})|(Y,Y_{j})_{{\bb L}_{2}(\omega)}|^{2}
<\infty\} \quad \mbox{ , with }\nonumber\\[-.4cm]
{~}\\[-.4cm]
(Y,Z)_{D({\bf B}^{\frac{m}{2}})}\,:=\,
\sum_{j=0}^{\infty}\,(1+\nu_{j}^{m})(Y,Y_{j})_{{\bb L}_{2}(\omega)}(Z,Y_{j})_{{\bb L}_{2}(\omega)}
\quad \,\forall \,Y,\,Z\,\in\,  D({\bf B}^{\frac{m}{2}})\nonumber
\end{eqnarray}
in accordance with {{Theorem B}} and  
{{Theorem C}} in the Appendix.\end{defi}
\begin{remark}\label{R5}   In the quotient space
  $D({\bf B}^{\frac{m}{2}}) / {\bb R}$ 
  one can also use the
scalar
product 
\begin{equation} \label{BeltrIISpaces} 
(Y,Z)_{*,{\frac{m}{2}}}\,:=\,
\sum_{j=1}^{\infty}\,(\nu_{j}^{m})(Y,Y_{j})_{{\bb L}_{2}(\omega)}(Z,Y_{j})_{{\bb L}_{2}(\omega)}
\quad \,\forall \,Y,\,Z\,\in\,  D({\bf B}^{\frac{m}{2}}) 
\,\,
\end{equation}
which is only a sesqui-scalar 
product on $D({\bf B}^{\frac{m}{2}})$.\end{remark}
\begin{remark}\label{R6}   The restriction of the Laplace-Beltrami operator ${\bf B}$ on quotient space
$D({\bf B}) / {\bb R}$ with the range ${\bb L}_{2}(\omega)/ {\bb R}$
is bijectiv. The inverse operator of this restriction is either 
explainable as the convolution operator with the Green function or
otherwise via Fourier series in surface spherical harmonics
(cf. {{Theorem C}} in the appendix).
\end{remark}
\begin{notationb}\label{N11}  To highlight the use of the restriction of  ${\bf B}$, we write\\[.1cm]
 ${\bf B}_{\bb R}\,:\,D({\bf B}) / {\bb R}\rightarrow {\bb L}_{2}(\omega)/ {\bb R}$ 
 if necessary. We will always denote the inverse of \,${\bf B}_{\bb R}$ by \,${\bf B}_{\bb R}^{-1}$.\end{notationb}
\begin{defi} \label{D5} We call ${\underline{\mathfrak{t}}}$ and ${\underline{\mathfrak{p}}}$
defined by: 
\begin{equation} 
\label{ToPo1}
{\underline{\mathfrak{t}}}\,:=\,curl({\psi{\;\!}\underline{x}})\,=\,grad{\;\!}{\psi}
{\;\!}\times{\;\!}  {{\underline{x}}}  \qquad\mbox{and}\qquad {\underline{\mathfrak{p}}}\,:=\,curl(curl({\chi{\;\!}\underline{x}}))
\end{equation}
toroidal and poloidal vector field.
The scalar functions ${\psi}$ and ${\chi}$ are denoted toroidal and
poloidal potentials, respectively.
\end{defi}
\begin{remark}\label{R7}  
Consider the special case 
${\underline{u}}({\underline{x}})\,:=\,{\underline{x}}$
and ${\underline{v}}({\underline{x}})\,:=\,{\eta}({\underline{x}}){\underline{x}}$.
Because of $curl{\;\!}{\underline{x}}\,:=\,{\underline{0}}$  it holds that 
${{\underline{x}}}=r{\underline{\mathfrak{e}}}_{r}$ is perpendicular to
$curl{\;\!}{\eta}{\underline{x}}\,=\,{\eta}{\;\!}{\;\!}curl{\;\!}
{\underline{x}}\,+grad{\;\!}{\eta} {\;\!}\times{\;\!} {{\underline{x}}}\,=\,
grad{\;\!}{\eta}{\;\!}\times{\;\!}  {{\underline{x}}}$.	
This justifies (\ref{ToPo1}).\end{remark}

\section{Decomposition\label{sec_dec}}
\subsection{Basic ideas of decomposition\label{sec_bas_dec}}
To achieve our goal to find a well-defined decomposition
of a divergence-free vector fields into toroidal and poloidal fields
in different Banach spaces, we have 
to ensure the uniqueness of toroidal and poloidal fields in these
spaces. Whereas we are going to use vanishing mean values on spherical surfaces
as the restriction for uniqueness, we will show the idea behind this
method in what follows.

We express this essential property in the characterization of solenoidal 
fields  on $\Omega_{\sigma}$, $0\leq{\sigma}<1$ as a Lemma,
where we formulate its statement in a slightly more general way.
\begin{lemmab} \label{LE1}
Let $\Omega$ be the open unit ball $\Omega$ in ${\bb R}^{n}$
for $n\,\in\,{\bb N}$, ${\underline{x}}\,\in\,{\bb R}^{n}$ 
and $r$ 
the
radius-coordinate in a
spherical coordinate system.
Then there exists no solenoidal vector function 
${\underline{g}} \in (C^{1}(\Omega))^{n}$, which is only depending of 
the variable $r$, with ${\underline{g}}(r)\,\neq\,{\underline{c}}$
$\forall \,r\,\in \, (0,1)$. Here ${\underline{c}}$ denotes any 
constant vector.
\end{lemmab} 
\begin{proof} The proof is by contradiction.
Firstly, the assertion of the Lemma is trivial for $n=1$.
So we assume without loss of generality that $ n\,\geq \,2$. 
Let the vector field ${\underline{g}}$
be written in Cartesian coordinates
\begin{equation*}
{\underline{g}}^{T}(r)\,=\,
(g_{1}(r),g_{2}(r),\dots,
g_{n}(r)) \,.
\end{equation*}
Since ${\underline{g}}(r)$ is in $(C^{1}(.))^{n}$, such that ${\underline{g}}(r)\,\neq\,{\underline{c}}$ there has to be 
at least one ${g}_{j^{*}}(r)\,,
j^{*}\,\in \{1,\dots,n \}$, with $\frac{dg_{j^{*}}(r)}{dr}\neq 0$ $(\circ)$.
In addition, the gradient of $r\,=\,r({\underline{x}})$
is a point on the unit sphere (for $r\neq 0$)
\begin{equation}
{\underline{\nabla}}_{\underline{x}}\,r\,=\,\frac{1}{r}{\underline{x}},\,
\,r\neq 0 \,\,,\nonumber
\end{equation}
what means that the $\frac{x_{j}}{r}$ are simply products of sinusoidal
functions of the angles
 $\vartheta_{1},\vartheta_{2},...,\varphi=\vartheta_{n-1}$.
Because of $(\circ)$ we get
$\frac{x_{j^{*}}}{r}\cdot\frac{dg_{j^{*}}(r)}{dr}\not\equiv 0$ and

\begin{equation}
div\, {\underline{g}} \,:=\,{\underline{\nabla}}_{\underline{x}}^{T}\cdot {\underline{g}}(r)\,
=\,\sum_{j=1}^{n}
\frac{x_{j}}{r}\cdot\frac{dg_{j}(r)}{dr}\,\not\equiv\,0
\quad
\nonumber
\end{equation}
which contradicts our assumption \, ${\textstyle{div}}\, {\underline{g}} \,=\,0$.
\end{proof}
\begin{remark}\label{R8}  Lemma \ref{LE1} is also true for all 
$\Omega_{\sigma}$,
since $\Omega_{\sigma}\,\subset\,\Omega_{o}$
for all ${\sigma}$
with  $0\,< \,\sigma \,\leq \,{\sigma}_{*} < 1$.
\end{remark}
\begin{remark}\label{R9}  The crucial point of the Lemma is that 
 in the 
Fourier series (developed in surface spherical 
harmonics) of solenoidal fields
on $\Omega_{\sigma}$
all terms only depending on the variable $r$ must 
be constants. This holds, e.g., for the coefficients of the first eigenfunction of the
Laplace-Beltrami operator $\bf B$.
\end{remark}
\noindent Now, we establish the equations for the decomposition of a
solenoidal vector field 
${\underline{u}}$ on ${\Omega}_{\sigma}$ into a toroidal field ${\underline{\mathfrak{t}}}$ and a poloidal field ${\underline{\mathfrak{p}}}$,\,
${\underline{u}}({\underline{x}})\,=\,{\underline{\mathfrak{t}}}
({\underline{x}})\,+\,
{\underline{\mathfrak{p}}}({\underline{x}})$.
Let the fields ${\underline{\mathfrak{t}}}$ and ${\underline{\mathfrak{p}}}$
be defined by (\ref{ToPo1}) in {{{Definition}}} \ref{D5},
i.e. by the conditions:
\begin{eqnarray} 
\label{ToPo2} div{\;\!}{\underline{\mathfrak{t}}}\,=\,0\,\,\mbox{and}\,\,
{{\underline{x}}}^{T}{\underline{\mathfrak{t}}}\,=\,0
 & \mbox{as well as} & div{\;\!}{\underline{\mathfrak{p}}}\,=\,0\,\,\mbox{and}\,\,
{{\underline{x}}}^{T}curl{\;\!}({\underline{\mathfrak{p}}})\,=\,0 \,\,.
\end{eqnarray}
We formulate the  conditions for uniqueness of ${\underline{\mathfrak{t}}}$ and  ${\underline{\mathfrak{p}}}$ by the restriction to vanishing mean values on the spherical surfaces $\omega_{r}$ for $r\,\in\,[\sigma,1],\, 0\,< \sigma<1\,$
(cf. {{Notation }}\ref{N1}):
\begin{equation}
\label{ToPo3}
\frac{1}{|\omega_{r}|}\intq_{\omega_{r}} 
{\underline{x}}^{T}{\underline{u}} d \omega_{r} \,=\,0\quad \mbox{and}
\quad 0\,= 
\frac{1}{|\omega_{r}|}\intq_{\omega_{r}} f \cdot {\psi}({\underline{x}})d \omega_{r}
\,=\,
\frac{1}{|\omega_{r}|}\intq_{\omega_{r}}  f \cdot {\chi}({\underline{x}})d \omega_{r}
\end{equation}
$\forall\,f\,\in\,{\bb C}[\sigma,1]$ with  $f\,=\,f(r)$. One can
understand the equations of ({\ref{ToPo3}}) as consequence 
of  {{Lemma}} \ref{LE2} in the next subsection.
\subsection{The Theorem of decomposition\label{sec_thm_dec}}
\begin{theorem} \label{Dec_Thm}
Let ${\underline{u}}\,\in\,{\underline{\bb W}}_{\,2}^{2}({\Omega}_{\sigma})$ be a solenoidal vector 
field defined on ${\Omega}_{\sigma}$ written in spherical coordinates, which
fulfills the first equations of ({\ref{ToPo3}}). Then there exist a toroidal vector field
${\underline{\mathfrak{t}}}$ and a poloidal vector field
${\underline{\mathfrak{p}}}$ which are uniquely determined by
({\ref{ToPo1}}), ({\ref{ToPo2}}), and ({\ref{ToPo3}}) with
${\underline{u}}\,=\,{\underline{\mathfrak{t}}}
\,+\,{\underline{\mathfrak{p}}}$ (at least in the sense of $
{\underline{\bb L}}_{2}({\Omega}_{\sigma})$).
The regularity of ${\underline{\mathfrak{t}}}$ and 
${\underline{\mathfrak{p}}}$ is inferred by that of ${\psi}$ and
${\chi}$.
\end{theorem}
\begin{remark}\label{R10}  The regularity of ${\psi}$ and ${\chi}$ is
inferred by the regularity of ${\underline{x}}^{T} curl({\underline{u}})$
and ${\underline{x}}^{T}{\underline{u}}$. There one has to study subspaces
of ${\underline{\bb L}}_{2}({\Omega}_{\sigma})$ as weighted Sobolev spaces
${\underline{\bb W}}_{2,r}^{m}(\sigma,1;D({\bf B}^{\frac{k}{2}}) / {\bb R})$
with $m,\,k\,\in\,{\bb N}$.\end{remark}
\begin{proof} 
(i) By simple calculations one finds ${\underline{\mathfrak{t}}}$ and
${\underline{\mathfrak{p}}}$ written in spherical 
coordinates as 
\begin{equation} \label{Toexpl}
{\underline{\mathfrak{t}}}_{\mathfrak{s}}\,=\&(0,\,
{\frac{\displaystyle{1}}{\displaystyle{\sin \vartheta}}}{\frac{\displaystyle{\partial {\psi}}}{\displaystyle{\partial {\varphi}}}},
\,-{\frac{\displaystyle{\partial {\psi}}}{\displaystyle{\partial {\vartheta}}}})^{T}_{\mathfrak{s}}\,\,,
\end{equation}
\begin{equation}\label{Poexpl}
{\underline{\mathfrak{p}}}_{\mathfrak{s}}\,=\,
{\frac{\displaystyle{1}}{\displaystyle{r}}}({\bf B}{\chi},\,
{\frac{\displaystyle{\partial}}{\displaystyle{\partial \vartheta}}}{\frac{\displaystyle{\partial (r {\chi})}}{\displaystyle{\partial {r}}}},\,{\frac{\displaystyle{1}}{\displaystyle{\sin \vartheta}}}
{\frac{\displaystyle{\partial}}{\displaystyle{\partial \varphi}}}{\frac{\displaystyle{\partial (r {\chi})}}{\displaystyle{\partial {r}}}}
\,)^{T}_{\mathfrak{s}}\,=\,
(-\-r\cdot\triangle{\;\!}{\chi},0,0)^{T}_{\mathfrak{s}}\,+\,
(grad{\;\!}({\frac{\displaystyle{\partial (r {\chi})}}{\displaystyle{\partial {r}}}}))_{\mathfrak{s}}\,\,.
\end{equation}
Taking the $curls$  of toroidal fields ${\underline{\mathfrak{t}}}$ and poloidal fields ${\underline{\mathfrak{p}}}$
is changing the types of the fields, i.e.
for more regular ${\underline{\mathfrak{t}}}\,,\,{\underline{\mathfrak{p}}}
\,\in\,{\underline{\bb C}}^{1}({\Omega}_{\sigma})$ the $curls$ of
toroidal fields are poloidal fields and the $curls$ of
poloidal fields are  toroidal field
\begin{equation}\label{curlToPo}
(curl{\;\!}{\underline{\mathfrak{t}}})_{\mathfrak{s}}
\,=\,
{\underline{{\tilde{\mathfrak{p}}}}}_{\mathfrak{s}}  \mbox{and} 
(curl{\;\!}{\underline{\mathfrak{p}}})_{\mathfrak{s}}
\,=\,
{\underline{{\tilde{\mathfrak{t}}}}}_{\mathfrak{s}}\,\,.
\end{equation}
Finally one obtains the following 
if one applies the 
$curl$ once more for ${\underline{\mathfrak{t}}}\,,\,{\underline{\mathfrak{p}}}
\,\in\,{\underline{\bb C}}^{2}({\Omega}_{\sigma})$ 
\begin{equation}\label{curl2ToPo}
(curl{\;\!}curl{\;\!}{\underline{\mathfrak{t}}})_{\mathfrak{s}}
\,=\,
{\underline{{\breve{\mathfrak{t}}}}}_{\mathfrak{s}}  \mbox{and} 
(curl{\;\!}curl{\;\!}{\underline{\mathfrak{p}}})_{\mathfrak{s}}
\,=\,
{\underline{{\breve{\mathfrak{p}}}}}_{\mathfrak{s}}\,\,.
\end{equation}
The relations ({\ref{curlToPo}}) and ({\ref{curl2ToPo}}) are also true 
for ${\underline{\mathfrak{t}}}\,,\,{\underline{\mathfrak{p}}}
\,\in\,{\underline{\bb W}}_{\,2}^{2}({\Omega}_{\sigma})$ by density
arguments.
\\[.3cm]
(ii) In the second step we prove orthogonality properties
of toroidal and poloidal fields in the ${\underline{\bb L}}_{2}({\Omega}_{\sigma})\,$- and 
${\underline{\bb W}}_{\,2}^{1}({\Omega}_{\sigma})\,$- sense. 
For arbitrary toroidal and poloidal fields
${\underline{\mathfrak{t}}}$ and ${\underline{\mathfrak{p}}}$ we see
\begin{eqnarray} \label{L2kugTP}
\intq_{\Omega_{\sigma}} {\underline{\mathfrak{t}}}_{\mathfrak{s}}^{T}\,
{\underline{\mathfrak{p}}}_{\mathfrak{s}} d {\underline{x}}& \,=\,&
\intq_{\Omega_{\sigma}} {\underline{\mathfrak{t}}}_{\mathfrak{s}}^{T}
(grad{\;\!}({\frac{\displaystyle{\partial (r {\chi})}}{\displaystyle{\partial {r}}}}))_{\mathfrak{s}}d {\underline{x}} \,=\,
\nonumber\\
& \,=\,&\intq_{\omega} {\underline{\mathfrak{t}}}_{\mathfrak{s}}^{T}
{\underline{x}}{\frac{\displaystyle{\partial (r {\chi})}}{\displaystyle{\partial {r}}}} d \omega  \,-\,{\sigma}^{-1}
\intq_{\omega_{\sigma}} {\underline{\mathfrak{t}}}_{\mathfrak{s}}^{T}{\underline{x}}
{\frac{\displaystyle{\partial (r {\chi})}}{\displaystyle{\partial {r}}}} d \omega_{\sigma}\,-\,
\intq_{\Omega_{\sigma}} div{\;\!}
{\underline{\mathfrak{t}}}_{\mathfrak{s}}
({\frac{\displaystyle{\partial (r {\chi})}}{\displaystyle{\partial {r}}}})
 d {\underline{x}} \,=\,\,0 \,,
\end{eqnarray}
where we have used partial integration and the conditions ({\ref{ToPo2}}) for ${\underline{\mathfrak{t}}}$.
In what follows it suffices to perform our investigations on scalar
functions and 2D-vector (tangential) fields both defined on the unit sphere $\omega$.
Here underlined small letters are to be 
understood as 2D-vector fields in the tangential plane.
Let us explain the Nabla operator ${\underline{\nabla}}_{{\;\!}\omega}$ on  $\omega$ by
\begin{equation}\label{nabsurf}
{\underline{\nabla}}_{{\;\!}\omega}\,:=\,
\,{\underline{\mathfrak{e}}}_{\vartheta}{\frac{\displaystyle{\partial {{~}}}}{\displaystyle{\partial {\vartheta}}}}
\,+\,{\underline{\mathfrak{e}}}_{\varphi}
{\frac{\displaystyle{1}}{\displaystyle{\sin \vartheta}}}{\frac{\displaystyle{\partial {{~}}}}{\displaystyle{\partial {\varphi}}}}\,\,\,.
\end{equation}
We define the surface-gradient, the surface-curl (also called surface-gradient-curl) and the
surface-divergence on $\omega$ for 
$\zeta\,\in\,{\bb C}^{1}({\omega})$, respectively, 
${\underline{b}}\,=\,(b_{\vartheta},b_{\varphi})^{T}\in\,{\underline{\bb C}}^{1}({\omega})$
through
\begin{eqnarray}\label{gradivsurf}
grad{\;\!}_{\omega}\zeta\,=\,{\underline{\nabla}}_{{\;\!}\omega}\zeta\,:=
\,{\underline{\mathfrak{e}}}_{\vartheta}{\frac{\displaystyle{\partial {\zeta}}}{\displaystyle{\partial {\vartheta}}}}
\,+\,{\underline{\mathfrak{e}}}_{\varphi}{\frac{\displaystyle{1}}{\displaystyle{\sin \vartheta}}}
{\frac{\displaystyle{\partial {\zeta}}}{\displaystyle{\partial {\varphi}}}}
&,&
{~}curl{\;\!}_{\omega}\zeta\,=\,{\underline{\nabla}}_{{\;\!}\omega} \,\times\,
 \zeta{\underline{\mathfrak{e}}}_{r}\,:=
\,{\underline{\mathfrak{e}}}_{\vartheta}{\frac{\displaystyle{1}}{\displaystyle{\sin \vartheta}}}
{\frac{\displaystyle{\partial {\zeta}}}{\displaystyle{\partial {\varphi}}}}
\,-\,{\underline{\mathfrak{e}}}_{\varphi}
{\frac{\displaystyle{\partial {\zeta}}}{\displaystyle{\partial {\vartheta}}}}{~}\,,
\nonumber\\[.1cm]
{~}&{~}&\hspace*{6cm}\mbox{ {~} }\\[-.4cm]
\mbox{ and }\hspace*{1cm}div{\;\!}_{\omega}{\underline{b}}\,=\,
{\underline{\nabla}}_{{\;\!}\omega}^{T}{\underline{b}}\,&:=&\,
{\frac{\displaystyle{1}}{\displaystyle{\sin \vartheta}}}(
{\frac{\displaystyle{\partial (b_{\vartheta}{\;\!} {\sin \vartheta})}}{\displaystyle{\partial {\vartheta}}}}\,+\,
{\frac{\displaystyle{\partial{\;\!}
                                                             b_{\varphi}}}{\displaystyle{\partial {\varphi}}}}) {~}. \nonumber 
\end{eqnarray}
Formally the Laplace-Beltrami operator ${\bf B}$ is
\begin{equation}\label{Beltrsurf}
{\bf B}\,=\,-{\underline{\nabla}}_{{\;\!}\omega}^{T}{\underline{\nabla}}_{{\;\!}\omega}
\qquad\mbox{ and }\qquad
{\bf B}{\;\!}{\zeta}\,=
\,-{\underline{\nabla}}_{{\;\!}\omega}^{T}{\underline{\nabla}}_{{\;\!}\omega}{\zeta}
\end{equation}
for more regular scalar functions ${\zeta}$.
It is convenient to use the Gauss theorem on $\omega$.
Similar to \cite{Chand} p.623-625 one finds 
\begin{equation}\label{Gauss_surf}
\intq_{\omega}{\;\!}div{\;\!}_{\omega}{\underline{b}}\,d \omega  \, 
=\,\int_{0}^{2 \pi}\int_{0}^{\pi}{\;\!}
({\frac{\displaystyle{\partial (b_{\vartheta}{\;\!} 
{\sin \vartheta})}}{\displaystyle{\partial {\vartheta}}}}
\,+\,
{\frac{\displaystyle{\partial{\;\!} b_{\varphi}}}{\displaystyle{\partial {\varphi}}}})
d {\vartheta} d {\varphi} \,=\,
2\pi[b_{\vartheta}{\;\!} {\sin \vartheta}]_{0}^{\pi}\,+\,\pi[b_{\varphi}]_{0}^{2\pi}
\,=\,0\,.
\end{equation}
The calculation rules for the Nabla-Calculus are collected in the
Appendix and hold true for the Nabla operator
${\underline{\nabla}}_{{\;\!}\omega}$ as well. In particular one verifies  
the property $(ii)$ in \eqref{nab_prod} as follows
\begin{equation}{\label{DivSuPro}}
div{\;\!}_{\omega}({\zeta}{\underline{b}}) = 
\,{\zeta} \,div{\;\!}_{\omega}{\underline{b}}
\,\,+\,\,{\underline{b}}^{T}\cdot\,grad{\;\!}_{\omega}{\zeta}\,\,\,\,. 
\end{equation}
Now, we provide an essential tool for the proof of
the orthogonality pro\-perties for our solenoidal fields.
In the following considerations it is advantageous to use the notation
of the system of
real-valued eigenfunctions of ${\bf B}$ (cf. {{Theorem B}} in the Appendix)
in the unified notation
\begin{equation}{\label{SuSphZ}}
\{Z_{l}^{-l},\dots,Z_{l}^{-1},Z_{l}^{0},Z_{l}^{1},\dots,Z_{l}^{l}\}
_{l=0}^{\infty}\,:=\,
\{{\tilde{Y}}_{l}^{l},\dots,{\tilde{Y}}_{l}^{1},
{{Y}}_{l}^{0},{{Y}}_{l}^{1},\dots,
{{Y}}_{l}^{l}\}
_{l=0}^{\infty}\,\,\,,
\end{equation}
where $\forall \, l=0,1,2,\dots\,$ the $\{Z_{l}^{k}\}_{k=-l}^{l}$ span the 
eigenspace to the eigenvalue $\nu_{(l)}:=l(l+1)$ in
the multiplicity
$N(\nu_{(l)}):=2l+1$.
In equation ({\ref{DivSuPro}}) we choose
${\underline{b}}\,=\,grad{\;\!}_{\omega}{Z_{l}^{k}}$ and
${\zeta}\,=\,{Z_{\tilde{l}}^{\tilde{k}}}$. Applying the Gauss theorem on
$\omega$ as in ({\ref{Gauss_surf}})  
we obtain
\begin{equation}\label{Gauss_toolA}
0\,=\,\intq_{\omega}{\;\!}div{\;\!}_{\omega}
({{Z_{\tilde{l}}^{\tilde{k}}}\,grad{\;\!}_{\omega}{Z_{l}^{k}}})\,d \omega  
 \,=\,\intq_{\omega}{\;\!}
(grad{\;\!}_{\omega}{Z_{\tilde{l}}^{\tilde{k}}})^{T}grad{\;\!}_{\omega}{Z_{l}^{k}}\,d \omega  \,
-\,\intq_{\omega}{\;\!}{Z_{\tilde{l}}^{\tilde{k}}}{\;\!}{\bf B}{Z_{l}^{k}}\,d \omega \,\,,
\end{equation}
and consequently,
\begin{equation}\label{Gauss_toolB} 
\intq_{\omega}{\;\!}
(curl{\;\!}_{\omega}{Z_{\tilde{l}}^{\tilde{k}}})^{T}curl{\;\!}_{\omega}{Z_{l}^{k}}\,d \omega 
 \,=\,
 \intq_{\omega}{\;\!}
(grad{\;\!}_{\omega}{Z_{\tilde{l}}^{\tilde{k}}})^{T}grad{\;\!}_{\omega}{Z_{l}^{k}}\,d \omega 
 \,=\,l(l+1)\delta_{l,{\tilde{l}}}{\;\!}\delta_{k,{\tilde{k}}}\quad,
\end{equation}
where we have used that the functions $\{\{Z_{l}^{k}\}_{k=-l}^{l}\}_{l=0}^{\infty}$
form a complete orthonormal system of Laplace-
Beltrami eigenfunctions 
in ${{\bb L}_{2}(\omega)}$ to the eigenvalues $\nu_{(l)}:=l(l+1)$ and
({\ref{gradivsurf}}) for $curl{\;\!}_{\omega}$.
\\[.3cm] 
Moreover, considering the ${{\underline{\bb L}}_{2}(\omega)}$ 
scalar product of surface-gradients and the surface-curls of the
$\{\{Z_{l}^{k}\}_{k=-l}^{l}\}_{l=0}^{\infty}$ one sees that these vector fields
are ${{\underline{\bb L}}_{2}(\omega)}$-orthogonal
\begin{eqnarray}\label{Gauss_toolC} 
\intq_{\omega}{\;\!}
(grad{\;\!}_{\omega}{Z_{l}^{k}})^{T}
curl{\;\!}_{\omega}{Z_{\tilde{l}}^{\tilde{k}}}\,d \omega 
 \,=\,\hspace*{10cm}\nonumber\\[-.3cm] 
\,=\,\intq_{\omega}{\;\!}{\big(}{\frac{\displaystyle{\partial  {Z_{l}^{k}}}}{\displaystyle{\partial {\vartheta}}}}
{\frac{\displaystyle{\partial {Z_{\tilde{l}}^{\tilde{k}}}}}{\displaystyle{\partial {\varphi}}}}\;-\;
{\frac{\displaystyle{\partial {Z_{l}^{k}}}}{\displaystyle{\partial {\varphi}}}}
{\frac{\displaystyle{\partial {Z_{\tilde{l}}^{\tilde{k}}}}}{\displaystyle{\partial {\vartheta}}}}{\big)}
{\frac{\displaystyle{d \omega}}{\displaystyle{\sin \vartheta}}}
 \,=\,
\intq_{\omega}{\;\!}{\big(}(-{\tilde{k}}){\frac{\displaystyle{\partial {Z_{l}^{k}}}}{\displaystyle{\partial {\vartheta}}}}
{Z_{\tilde{l}}^{- \tilde{k}}}\;+\;k
{Z_{l}^{- k}}
{\frac{\displaystyle{\partial {Z_{\tilde{l}}^{\tilde{k}}}}}{\displaystyle{\partial {\vartheta}}}}{\big)}
{\frac{\displaystyle{d \omega}}{\displaystyle{\sin \vartheta}}}\,=\hspace*{.5cm}\\
\hspace*{-1.5cm}\,=\,\pi |k|(\sqrt{2})^{-\delta_{k,0}}
\sqrt{\frac{\displaystyle{(l-|k|)!(2l+1)}}{\displaystyle{(l+|k|)!2\pi}}}
\sqrt{\frac{\displaystyle{({\tilde{l}}-|k|)!(2{\tilde{l}}+1)}}
{\displaystyle{({\tilde{l}}+|k|)!2\pi}}}{\delta}_{k,-{\tilde{k}}}
{\big[}(P_{l}^{|k|}P_{{\tilde{l}}}^{|k|})(\cos {\vartheta}){\big]_{0}^{\pi}} 
\,=\,0\,,\nonumber
\end{eqnarray}
where one shows the vanishing of the integrals by considering that
$P_{l}^{|k|}(\pm 1)\:\!=\:\!P_{{\tilde{l}}}^{|k|}(\pm 1)\:\!=\:\!0$ 
for $0<|k|\leq
\min(l,{\tilde{l}})$.
Let the toroidal and poloidal scalar potentials ,${\psi}$ respectively ${\chi}$, be formally expanded in 
partial Fourier's series in surface spherical harmonics as
\begin{eqnarray}\label{Expkug}
{\psi}({\underline{x}})\,=\,\sum _{l=1}^{\infty}\sum _{k=-l}^{l}{\;\!}{\hat{\psi}}_{l,k}(r)Z_{l}^{k}(\vartheta,\varphi)
 & \mbox{and} &
{\chi}({\underline{x}})\,=\,\sum _{l=1}^{\infty}\sum_{k=-l}^{l}{\;\!}{\hat{\chi}}_{l,k}(r)Z_{l}^{k}(\vartheta,\varphi)
\,\,,
\end{eqnarray}
where we  have used the restrictions ({\ref{ToPo3}}) by the prescribed  
${\hat{\psi}}_{0,0}(r)\,=\,{\hat{\chi}}_{0,0}(r)\,:=\,0$. 
The series in ({\ref{Expkug}}) make already sense for allmost all
$r\,\in\,(\sigma,1)$ if one supposes 
${\psi}\,,\,{\chi}\,\in{\underline{\bb L}}_{2}({\Omega}_{\sigma})$.
\\[.3cm]
For the moment we ignore questions about the convergence of the series
and the regularity of 
their coefficients.
Employing the elements of these series we explain  $\forall \; l\,\in\,{\bb N}$ and $\forall \;
k\,=-l,\dots,l$
\begin{eqnarray} 
\label{lkToPo}
{\underline{\mathfrak{t}}}^{(l,k)}_{\mathfrak{s}}({\underline{x}})\,:=\,(curl({\hat{\psi}}_{l,k}(r)Z_{l}^{k}(\vartheta,\varphi){\;\!}
))_{\mathfrak{s}} & \mbox{and} &
{\underline{\mathfrak{p}}}^{(l,k)}_{\mathfrak{s}}({\underline{x}}) \,:=\,(curl(curl({\hat{\chi}}_{l,k}(r)Z_{l}^{k}(\vartheta,\varphi){\;\!}
)))_{\mathfrak{s}}\,
\end{eqnarray}
as  purpose-built toroidal fields and poloidal fields.
In what follows we suppose, that the sequences of coefficients
\begin{eqnarray} 
\label{ColkToPo}
\{\{{\hat{\psi}}_{l,k}(r)\}_{k=-l}^{l}\}_{l=1}^{\infty} & \mbox{and}
 &\{\{{\hat{\chi}}_{l,k}(r)\}_{k=-l}^{l}\}_{l=1}^{\infty}\,\,\,
\end{eqnarray}
are elementwise arbitrarily defined $\forall\,r\,\in [\sigma,1]$ and
$\forall\,r\,\in[0,1]$, respectively.
In view of Fubini's theorem we are going to show, that all the vector
fields in the general system
of solenoidal vector functions
\begin{eqnarray} 
\label{lkToPo}
\{\{{\underline{\mathfrak{t}}}^{(l,k)}_{\mathfrak{s}}\}_{k=-l}^{l}\}_{l=1}^{\infty}
\cup \{\{{\underline{\mathfrak{p}}}^{(l,k)}_{\mathfrak{s}}\}_{k=-l}^{l}\}_{l=1}^{\infty}
\end{eqnarray}\\[-.3cm]
$\forall\,r\,\neq \,0$  already form an orthogonal family of vector functions in the sense of ${\underline{\bb L}}_{2}({\omega})$.
At first we check, that the toroidal fields in ({\ref{lkToPo}}) are
mutually orthogonal:
\begin{eqnarray}\label{Orth_T_surf}
\intq_{\omega}({\underline{\mathfrak{t}}}^{({\tilde{l}},{\tilde{k}})}_{\mathfrak{s}})^{T}
{\underline{\mathfrak{t}}}^{(l,k)}_{\mathfrak{s}}\,d \omega  \, 
=\,{\hat{\psi}}_{{\tilde{l}},{\tilde{k}}}(r){\hat{\psi}}_{l,k}(r)\intq_{\omega}{\;\!}
(curl{\;\!}_{\omega}{Z_{\tilde{l}}^{\tilde{k}}})^{T}
 curl{\;\!}_{\omega}{Z_{l}^{k}}\,d\omega  \,
\,=\,l(l+1){\hat{\psi}}_{{\tilde{l}},{\tilde{k}}}{\hat{\psi}}_{l,k}\delta_{l,{\tilde{l}}}
{\;\!}\delta_{k,{\tilde{k}}} \,\, 
\end{eqnarray}
$\forall \; l,{\tilde{l}}\,\in\,{\bb N}$ and $\forall \;
k,{\tilde{k}}\,=-l,\dots,l;\,{\tilde{k}}\,=-{\tilde{l}},\dots,{\tilde{l}}$, where we have used ({\ref{Toexpl}}). 
Similarly one sees for the poloidal fields in ({\ref{lkToPo}}) that 
\begin{eqnarray} \label{Orth_P_surf}
r^{2}\intq_{\omega}{\;\!}
({\underline{\mathfrak{p}}}^{({\tilde{l}},{\tilde{k}})}_{\mathfrak{s}})^{T}
{\underline{\mathfrak{p}}}^{(l,k)}_{\mathfrak{s}}\,d \omega  \, 
&=&\,{\hat{\chi}}_{{\tilde{l}},{\tilde{k}}}(r){\hat{\chi}}_{l,k}(r)
\intq_{\omega}{\;\!}{{\bf B}}{Z_{\tilde{l}}^{\tilde{k}}}
{{\bf B}}{Z_{l}^{k}}\,d \omega \,\nonumber\\
{~}&{~}&+
{\frac{\displaystyle{\partial (r{{\hat{\chi}}_{{\tilde{l}},{\tilde{k}}}(r)})}}{\displaystyle{\partial {r}}}}
{\frac{\displaystyle{\partial (r{{\hat{\chi}}_{l,k}(r)})}}{\displaystyle{\partial {r}}}}
\intq_{\omega}{\;\!}
(grad{\;\!}_{\omega}{Z_{\tilde{l}}^{\tilde{k}}})^{T}
 grad{\;\!}_{\omega}{Z_{l}^{k}}\, d \omega  \,=\,\,\,\\[-.2cm]
&=&\,l(l+1){\Big[}l(l+1){\hat{\chi}}_{{\tilde{l}},{\tilde{k}}}(r){\hat{\chi}}_{l,k}(r)
\,+\,{\frac{\displaystyle{\partial (r{{\hat{\chi}}_{{\tilde{l}},{\tilde{k}}}(r)})}}{\displaystyle{\partial {r}}}}
{\frac{\displaystyle{\partial (r{{\hat{\chi}}_{l,k}(r)})}}{\displaystyle{\partial {r}}}}
{\Big]}\delta_{l,{\tilde{l}}}{\;\!}\delta_{k,{\tilde{k}}}\,.\nonumber 
\end{eqnarray}
Analogously to ({\ref{L2kugTP}}) for arbitrary toroidal and poloidal
fields in ({\ref{lkToPo}}) one obtains 
that these fields are already orthogonal in the sense of
${\underline{\bb L}}_{2}({\omega})$, namely,
\begin{eqnarray} \label{Orth_T+P_surf}
r\intq_{\omega}{\;\!}({\underline{\mathfrak{p}}}^{(l,k)}_{\mathfrak{s}})^{T}
{\underline{\mathfrak{t}}}^{({\tilde{l}},{\tilde{k}})}_{\mathfrak{s}}
\,d \omega  \, 
& = &{\frac{\displaystyle{\partial (r{{\hat{\chi}}_{l,k}(r)})}}{\displaystyle{\partial {r}}}}
{\hat{\psi}}_{{\tilde{l}},{\tilde{k}}}(r)
\intq_{\omega}{\;\!}
(grad{\;\!}_{\omega}{Z_{l}^{k}})^{T}
curl{\;\!}_{\omega}{Z_{\tilde{l}}^{\tilde{k}}}\,d \omega \,\,=\,0\,,
\end{eqnarray}
where we have used ({\ref{Gauss_toolC}}).
Again, these relations are true $\forall \; l,{\tilde{l}}\,\in\,{\bb N}$ and $\forall \;
k,{\tilde{k}}\,=-l,\dots,l;\,{\tilde{k}}\,=-{\tilde{l}},\dots,{\tilde{l}}$. In
deriving ({\ref{Orth_P_surf}}) we have applied 
the explicit shape of the poloidal fields ({\ref{Poexpl}}), the tool
({\ref{Gauss_toolB}}) and the known eigenvalues $\nu_{(l)}:=l(l+1)$ of
the Operator ${{\bf B}}$.
In the next step we will investigate the System (\ref{lkToPo}) with
regard to the Dirichlet (scalar) product 
$(.,.)_{D}$ (cf. (\ref{Dirichlet})).
Take ${\underline{\mathfrak{t}}}^{(l,k)}_{\mathfrak{s}}$  and
${\underline{\mathfrak{p}}}^{(l,k)}_{\mathfrak{s}}$  as in
(\ref{lkToPo}). We will show that
$({\underline{\mathfrak{p}}}^{(l,k)}_{\mathfrak{s}},{\underline{\mathfrak{t}}}^{(l,k)}_{\mathfrak{s}})_{D}\,=\,0$. First 
we use partial integration to arrive at
\begin{eqnarray} \label{TPDirichlet} 
({\underline{\mathfrak{p}}}^{(l,k)}_{\mathfrak{s}},
{\underline{\mathfrak{t}}}^{({\tilde{l}},{\tilde{k}})}_{\mathfrak{s}})_{D} & = & -\,
\intq_{\Omega_{\sigma}}{\;\!}({\underline{\mathfrak{p}}}^{(l,k)}_{\mathfrak{s}})^{T}(
\triangle{\;\!}{\underline{\mathfrak{t}}}^{({\tilde{l}},{\tilde{k}})}_{\mathfrak{s}})
\,d  {\Omega_{\sigma}} \, +\\{~} &{~} &
+\,{\underbrace{\intq_{\omega}{\;\!}({\underline{\mathfrak{p}}}^{(l,k)}_{\mathfrak{s}})^{T}
({\underline{\nabla}}\cdot({\underline{\mathfrak{t}}}^{({\tilde{l}},{\tilde{k}})}_{\mathfrak{s}})^{T})^{T}
{\underline{\mathfrak{e}}}_{r}
\,d \omega}_{\displaystyle{=0 (*)}}} \,-\,{\underbrace{
\intq_{\omega_{\sigma}}{\;\!}({\underline{\mathfrak{p}}}^{(l,k)}_{\mathfrak{s}})^{T}
({\underline{\nabla}}\cdot({\underline{\mathfrak{t}}}^{({\tilde{l}},{\tilde{k}})}_{\mathfrak{s}})^{T})^{T}{\sigma}^{-1}
{\underline{\mathfrak{e}}}_{r}
\,d \omega_{\sigma}}_{\displaystyle{=0 (*)}}} \,.\nonumber
\end{eqnarray} 
Here
$({\underline{\nabla}}\cdot({\underline{\mathfrak{t}}}^{({\tilde{l}},{\tilde{k}})}_{\mathfrak{s}})^{T})^{T}$
denotes 
the Jacobian  matrix (in spherical polar coordinates) and ${\underline{\mathfrak{e}}}_{r}$ the unit vector in
${r}$-direction, where we use it in the sense of
${\underline{\mathfrak{e}}}_{r}\,=\,(1,0,0)_{\mathfrak{s}}^{T} 
$ here. 
With the relations (\ref{Toexpl}),  ( \ref{Poexpl}), \eqref{EF_Ver_Lap2}
and the representation
 ${\displaystyle{\triangle{\;\!}{\underline{\mathfrak{t}}}
^{({\tilde{l}},{\tilde{k}})}_{\mathfrak{s}}\,=\,(0,(\triangle{\;\!}{{\mathfrak{t}}}^{({\tilde{l}},{\tilde{k}})})_{\vartheta}
,(\triangle{\;\!}{{\mathfrak{t}}}^{({\tilde{l}},{\tilde{k}})})_{\varphi}
)^{T}_{\mathfrak{s}}}}$ we derive from (\ref{TPDirichlet})
\begin{eqnarray} \label{TPDirichlet1} 
({\underline{\mathfrak{p}}}^{(l,k)}_{\mathfrak{s}},
{\underline{\mathfrak{t}}}^{({\tilde{l}},{\tilde{k}})}_{\mathfrak{s}})_{D} & = &
\intq_{\Omega_{\sigma}}{\;\!}\{
{\overbrace{(r\cdot\triangle{\;\!}{\hat{\chi}}_{l,k}(r)Z_{l}^{k},0,0)_{\mathfrak{s}}\cdot
(\triangle{\;\!}{\underline{\mathfrak{t}}}^{({\tilde{l}},{\tilde{k}})}_{\mathfrak{s}})}^{\displaystyle{=0}}}\,-\,
((grad{\;\!}({\frac{\displaystyle{\partial (r {{\hat{\chi}}_{l,k}(r)Z_{l}^{k}})}}
{\displaystyle{\partial {r}}}}))_{\mathfrak{s}})^{T}(
\triangle{\;\!}{\underline{\mathfrak{t}}}^{({\tilde{l}},{\tilde{k}})}_{\mathfrak{s}})
\}\,d  {\Omega_{\sigma}} \nonumber\\
{~} & = &
\intq_{\Omega_{\sigma}}{\;\!} {\frac{\displaystyle{\partial (r {{\hat{\chi}}_{l,k}(r)Z_{l}^{k}})}}
{\displaystyle{\partial {r}}}} div{\;\!}(
\triangle{\;\!}{\underline{\mathfrak{t}}}^{({\tilde{l}},{\tilde{k}})}_{\mathfrak{s}})\,d  {\Omega_{\sigma}}
 -\intq_{\omega}{\;\!}{\underbrace{(
\triangle{\;\!}{\underline{\mathfrak{t}}}^{({\tilde{l}},{\tilde{k}})}_{\mathfrak{s}})^{T}
({\frac{\displaystyle{\partial (r {{\hat{\chi}}_{l,k}(r)Z_{l}^{k}})}}
{\displaystyle{\partial {r}}}})
{\underline{\mathfrak{e}}}_{r}}_{\displaystyle{=0}}}
\,d \omega\\
{~} & {~} & +\,\intq_{\omega_{\sigma}}{\;\!}{\underbrace{(
\triangle{\;\!}{\underline{\mathfrak{t}}}^{({\tilde{l}},{\tilde{k}})}_{\mathfrak{s}})^{T}
({\frac{\displaystyle{\partial (r {{\hat{\chi}}_{l,k}(r)Z_{l}^{k}})}}
{\displaystyle{\partial {r}}}}){\sigma}^{-1}
{\underline{\mathfrak{e}}}_{r}}_{\displaystyle{=0}}}
\,d \omega_{\sigma}\,.\nonumber
\end{eqnarray} 
(iii) It remains to prove that the equations for 
of ${\psi}$ and  ${\chi}$ have unique solutions, i.e.
\begin{eqnarray} 
\label{BewToPo} 
{\bf B}{\psi}\,=\,
{{\underline{x}}}^{T}(curl{\;\!}{\underline{u}})
 & \mbox{and} & {\bf B}{\chi}\,=\,
{{\underline{x}}}^{T}{\underline{u}}\,.\,\quad\hspace*{3cm} .
\end{eqnarray}\end{proof}
\noindent It is worth to  consider the following property of gradients of
 harmonic functions.
\begin{lemmab} \label{LE2} Let $p\, \in\,H_{\sigma}$ be a harmonic function on  $\Omega_{\sigma}$, $0\leq{\sigma}<1$.
Then the gradients ${\underline{\nabla}} \hspace*{.1cm} p$ of $p\,
\in\, H_{\sigma}$ are poloidal fields ${\underline{\mathfrak{p}}}$.
(The systems $\,H_{\sigma}$ of harmonic functions are explained in the
appendix.)\end{lemmab} 
\begin{proof}	The statement of the Lemma is a simple conclusion  of the application of (\ref{Poexpl}) on ${\chi}({\underline{x}})\,:=\,
p({\underline{x}})\,$, \\
where one uses:  $0\,=\,{\triangle}{\chi}\,=\,{\triangle}{p} \,\, \forall \,p\,
\in\,H_{\sigma}$  . \hfill \end{proof}
\begin{remark}\label{R_harmspan}  The statement of Lemma \ref{LE2} is also true for any $p\, \in\,{\mathrm{span}\,} H_{\sigma}$.
\end{remark}

\section{The eigenvalue problem for the Stokes operator\label{sec_sto}}
\subsection{Preliminary theoretical results\label{sec_theosto}}
In Subsection \ref{sec_spac_not} we have introduced the Stokes operator ${\bf A}_{\sigma}\,:=\,- \Upsilon \triangle $
as a linear operator from 
${\underline{\bb S}}^{2}({\Omega}_{\sigma})$ onto 
${\underline{\bb S}}({\Omega}_{\sigma})\,
\subset\,{\underline{\bb L}}_{2}({\Omega}_{\sigma})$. From another point of view,
one can understand the Stokes operator ${\bf A}_{\sigma}$ 
in the sense of Friedrichs' extension
using the restriction of the the Leray-Helmholtz projector $\Upsilon$ on 
${\underline{{\cal C}}}^{\infty}_{o}({\Omega_{\sigma}})$ 
 \begin{eqnarray}\label{DefLerHe_infty}
\Upsilon_{{\underline{{\cal C}}}^{\infty}_{o}
({\Omega_{\sigma}})}({\underline{v}}({\underline{x}}))\,:=\,
\frac{1}{4\pi}curl{\;\!}\intq_{{\Omega}_{\sigma}}{\;\!}
\frac{curl{\;\!}{\underline{v}}({\underline{y}})}
{||{\underline{x}}-{\underline{y}}||_{\cal{E}}}\,d {\Omega}_{\sigma}\,-\,
grad{\;\!}p({\underline{x}})\,\,\mbox{, with:} \,\,-\triangle{\;\!}p\,=\,0\,\,,
\end{eqnarray}
and multiply it with the negative Laplacian
on the domain: ${\underline{\cal V}}_{\sigma}\,:=\,  
\{{\underline{v}}\,\in\,
{\underline{{\cal C}}}^{\infty}_{o}({\Omega_{\sigma}}):div{\;\!}{\underline{v}}=0\}$
(cf. {{Notation \ref{N10}}}).
Using the well-known results on Friedrichs' extensions of symmetric, linear and semi-bounded operators, one obtains that the Stokes operator ${\bf A}_{\sigma}$
is a self-adjoint
positive definite operator with the energy space 
${\underline{\bb S}}^{1}({\Omega}_{\sigma})$.
Since the embedding of the energy space 
${\underline{\bb S}}^{1}({\Omega}_{\sigma})$ into the Hilbert space 
${\underline{\bb S}}({\Omega}_{\sigma})$ is compact, Rellich's Theorem
for operators with pure
point spectrum is applicable to the Stokes operator ${\bf A}_{\sigma}$, cf. \cite{Triebel} p.258. In particular the spectrum of every operator ${\bf A}_{\sigma}$ consists
only of eigenvalues with finite multiplicity. 
This property of ${\bf A}_{\sigma}$ is the precondition of the Theorem
on completeness (and existence), cf. \cite{Triebel} p.254. 
Let us collect the properties of the Stokes operator ${\bf A}_{\sigma}$.
\begin{theorem} \label{thmstok}
The Stokes operator ${\bf A}_{\sigma}$ is positive and self-adjoint.
Its inverse ${\bf A}_{\sigma}^{-1}$ is injective, self-adjoint and compact.
\end{theorem}
\noindent The proof of Theorem \ref{thmstok} 
is a simple modification of 
Theorems 4.3 and 4.4 in \cite{CoFoi}. The essential tools are
the Rellich theorem and Lax-Milgram lemma (for the existence and uniqueness of the weak solutions).
The well-known theorem of Hilbert and
regularity results like \cite[Prop.~I.2.2]{Temam}
lead to more detailed results.  
\begin{lemmab} \label{STOeiFU}
The Stokes operator  ${\bf A}_{\sigma}$ for any fixed $ {\sigma}:0\leq{\sigma}<1$ is an operator with a pure point spectrum.
The eigenvalues of the Stokes operator ${\bf A}_{\sigma}$: $\{\lambda_{j}\}_{j=1}^{\infty}$ are positive
and of finite multiplicity.
They can be ordered by their values, taking into account their multiplicities:$\{\lambda_{j}\}_{j=1}^{\infty}$, with 
${\lim}_{j \to \infty}\lambda_{j}\,=\,\infty$.
The associated eigenfunctions  $\{{\underline{w}}_{j}({\underline{x}})\}_{j=1}^{\infty}$ of
${\bf A}_{\sigma}$ are  complete in the range $R({\bf A}_{\sigma})\,=\,{\underline{\bb S}}({\Omega}_{\sigma})$ and 
in
$D({\bf A}_{\sigma})\,=\,\{{\underline{v}}\,\in\,{\underline{\bb S}}({\Omega}_{\sigma}):\,
\sum_{j=1}^{\infty}\,
\lambda_{j}^{2}|({\underline{v}},{\underline{w}}_{j})_{{\underline{\bb S}}(.)}|^{2}
<\infty\}\,=\,
{\underline{\bb S}}^{2}({\Omega}_{\sigma})$.
The eigenfunctions 
$\{{\underline{w}}_{j}({\underline{x}})\}_{j=1}^{\infty}$
of ${\bf A}_{\sigma}$ (counted in multiplicity)
are an orthogonal basis of
${\underline{\bb S}}(.)$ and ${\underline{\bb S}}^1(.)$ as well.
We obtain that
\begin{align*}
{{(a)}}&\quad  {\bf A}_{\sigma}{\underline{w}}_{j}=\lambda_{j}{\underline{w}}_{j}\quad\mbox{for
}\quad{\underline{w}}_{j}\in D({\bf A}_{\sigma})
\quad\forall\,j=1,2,\dots\,;\\
{(b)}& \quad
0\,<\lambda_{1}\leq\,\lambda_{2}\,\leq\cdots\leq\,\lambda_{j}
\,\leq\cdots\quad\mbox{and}\quad 
\lim_{j\rightarrow\infty}\lambda_{j}=\infty\,;
\\[-2mm]
{(c)}&\quad
\|{\underline{w}}_{j}\|_{{\underline{\bb 
S}}}\,=1\quad\forall \, j=1,2,\dots\,.
\end{align*}
\end{lemmab}	
\noindent Keep in mind that the following property for self-adjoint
operators with a pure point spectrum holds:
\begin{remark}\label{ROrtHilbert}  For any fixed $ 0\leq{\sigma}<1$
  the Stokes operator  ${\bf A}_{\sigma}$ is a  self-adjoint operator
  with pure point spectrum. So we can conclude that the Stokes
  eigenfunctions ${\underline{w}}_{j}$ and 
${\underline{w}}_{k}$ associated to the eigenvalues $\lambda_{j}$ and $\lambda_{k}$ with
$\lambda_{j}\,\neq\,\lambda_{k}$ are orthogonal in ${\underline{\bb S}}$.
\end{remark}
\subsection{The governing equations for the scalar potentials of the Stokes eigenfunctions \label{sec_sto_ew1}}
In this paragraph we are going to split the Stokes eigenvalue Problem
\ref{Prob1} into separate 
eigenvalue problems. More precisely, from the eigenvalue problems for
${\underline{u}}$, $\lambda$ and 
$p$ we derive boundary value problems
for the scalar potential $\psi$ of a toroidal eigenfunction
${\underline{u}}\,=\,{\underline{\mathfrak{t}}}\,=\,curl{\;\!}({\psi{\;\!}\underline{x}})$
and for the scalar potential $\chi$ of the correspondent poloidal
Stokes eigenfunction
${\underline{u}}\,=\,{\underline{\mathfrak{p}}}\,=\,curl{\;\!}curl{\;\!}({\chi{\;\!}\underline{x}})$. 
For fixed $0\leq{\sigma}<1$ we are looking for
solutions ${\underline{u}}$, $\lambda$ and $p$ 
of:\\[.2cm]
$- \triangle {{\;\!}}{\underline{u}}\,
 +\,{\underline{\nabla}} \hspace*{.1cm} p
 = \,\lambda {\underline{u}}\,\quad$ and
$\,\quad div \hspace*{.1cm}{\underline{u}}\,=\,
{\underline{\nabla}}^{T}\cdot {\underline{u}}\,=\,0\quad$
{in} ${\Omega}_{\sigma}$ , where
${\underline{u}}\,=\,{\underline{0}}$ on $\partial{\Omega}_{\sigma}$.\\[.2cm]
We note, that any eigenfunction ${\underline{u}}$ has to be a solenoidal vector field as a zero divergence field with
prescribed homogeneous Dirichlet boundary conditions on
$\partial{\Omega}_{\sigma}$ (cf. Notation \ref{N1}).
Additionally we observe, that the $curl$ of ${\underline{u}}$ is also a solenoidal 
vector field. This follows directly from
$div{\;\!}curl{\;\!}{\underline{u}}\,=\,0$ and since for $r>0$ the
surface $
\,{\omega}_{r}\,\subset \,{\overline{\Omega}}_{o}$ is closed
with empty boundary (cf. the Stokes theorem). 
Unless otherwise noted, we suppose whenever necessary $r>0$ and we assume higher regularity of $\psi{\;\!}$ as well as of $\chi{\;\!}$
in the following.\\[.2cm]
In the first step we derive the boundary value problem for the scalar
potential $\psi$ of a toroidal
eigenfunction. 
For this 
we apply the $curl$ to the equation $- \triangle {{\;\!}}{\underline{u}}
 \,+ \,{\underline{\nabla}} \hspace*{.1cm} p
\, = \,\lambda {\underline{u}}$ to obtain
 $-\,curl{\;\!}\left(\triangle{\;\!}{\underline{u}}\right) =
-\triangle{\;\!}\left(curl{\;\!}{\underline{u}}\right)=
\lambda{\;\!}(curl{\;\!} 
{\underline{u}})$ \,.
By a simple calculation for the radial component of $curl{\;\!}
{\underline{u}}$ we find 
\begin{eqnarray} 
\label{EF__Equ_Tor1_Rad} 
\hspace*{3cm}
- \triangle{\;\!}\left({\underline{x}}^{T} (curl{\;\!} 
{\underline{u}})\right)
&\,=\,&\lambda{\;\!}\left({\underline{x}}^{T} (curl{\;\!} 
{\underline{u}})\right)\quad .
\end{eqnarray} 
Using ${\bf B}{\psi}\,=\,
{{\underline{x}}}^{T}(curl{\;\!}{\underline{u}})$ (equation
(\ref{BewToPo})) we obtain the 
equation for the scalar function  ${\psi}$ as
\begin{eqnarray} \label{EF__Equ_Tor1_Rad1} {~}\hspace*{3.8cm}
- {\bf B}(\triangle{\;\!}{\psi})
&\,=\,& {\bf B}(\lambda{\;\!}{\psi})\quad ,
\end{eqnarray} 
where we have used (\ref{Lapl_Komm-Belt}).  Finally we 
apply ${\bf B}_{\bb R}^{-1}$ to  (\ref{EF__Equ_Tor1_Rad1}) and find
\begin{eqnarray} \label{EF__Equ_Tor1_Rad2} 
{~}\hspace*{3cm}\quad\quad\quad\quad
- \triangle{\;\!}{\psi}
&\,=\,& \lambda{\;\!}{\psi}\quad ,
\end{eqnarray} \\[-.4cm]
since for all fixed $ r> 0\:,r\in  [{\sigma},1]$ taken as a parameter
one can assure the property ${\psi}(r,.,.)\,\in\,D({\bf B}_{\bb
  R}^{-1})\,=\,{\bb L}_{2}(\omega)/ {\bb R}$ (cf. Notation \ref{N11})
and (\ref{ToPo3})). 
The validity of (\ref{EF__Equ_Tor1_Rad2}) for $r=0$  is ensured  by
density arguments.
\\[.3cm]
We observe that the boundary conditions for ${\psi}$ are imposed  by the prescribed homogeneous 
Dirichlet boundary conditions for a toroidal eigenfunction ${\underline{u}}\,=\,{\underline{\mathfrak{t}}}$ on $\partial{\Omega}_{\sigma}$. 
To gather ${\psi}\,=\,0 $ on $\partial{\Omega}_{\sigma}$
from ${\underline{u}}\,=\,grad{\;\!}{\psi}
\,\times\, {{\underline{x}}}\,=\,{\underline{0}}$ on $\partial{\Omega}_{\sigma}$ (cf. (\ref{ToPo1})) 
(this is equivalent to
${\psi}\,=\,const. $ on $\partial{\Omega}_{\sigma}$ , cf. \cite{GirRav} p. 35, Remark 2.3)
one has to use the condition  (cf. (\ref{ToPo3})) \\[.2cm]
\hspace*{1.3cm} $\frac{1}{|\omega_{r}|}\intq_{\omega_{r}} {\psi}({\underline{x}})d \omega_{r}\,=\,0$ 
\quad\quad for
$r\,=\,1$
and $r\,=\,\sigma > 0$ iff  ${\Omega}_{\sigma}$ is a spherical shell.\\[.1cm]
We summarize the eigenvalue problem for the scalar potential $\psi$ of the correspondent toroidal Stokes eigenfunction ${\underline{u}}\,=\,curl({\psi{\;\!}\underline{x}})$
in the following problem, where the eigenvalues $\lambda$ coincide.
\begin{problem} \label{Prob3} For $0\leq{\sigma}<1$ we seek solutions
  $\psi$  and  $\lambda$  
 fulfilling\\[.3cm]
$- \triangle {{\;\!}}{\psi}\,
 = \,\lambda {\psi}\,\, $ {in} ${\Omega}_{\sigma}$ , with 
$\psi\,=\,0$ on $\partial{\Omega}_{\sigma}$ and 
$\frac{1}{|\omega_{r}|}\intq_{\omega_{r}} f \cdot {\psi}({\underline{x}})d \omega_{r}\,=\,0$
$\forall r >0:r\,\in\,[\sigma,1]$ (cf. (\ref{ToPo3})).\\[.05cm] 
The last condition has to be satisfied $\forall \,f\,=\,f(r)$. 
\end{problem}
\noindent Secondly from the Stokes eigenvalue problem for
a poloidal eigenfunction 
${\underline{u}}\,=\,{\underline{\mathfrak{p}}}$ we construct the
eigenvalue problem for the 
corresponding scalar poloidal potential $\chi$ in the shape of a boundary value 
problem with the Bi-Laplacian, where the eigenvalues $\lambda$ will
coincide again. Similarly to the toroidal eigenfunction we obtain via
the 
application of the $curl$ on the equation\\[.2cm]
 $-\,curl{\;\!}\left(\triangle{\;\!}{\underline{u}}\right) =
-\triangle{\;\!}\left(curl{\;\!}{\underline{u}}\right)=
\lambda{\;\!}(curl{\;\!} 
{\underline{u}})$ \quad \quad: \quad $-\,curl{\;\!} curl{\;\!}\left(\triangle{\;\!}{\underline{u}}\right) =
-\triangle{\;\!}\left(\triangle{\;\!}{\underline{u}}\right)=
\lambda{\;\!}(\triangle{\;\!}
{\underline{u}})$.\\[.2cm] 
By an simple calculation for the radial component of ${\underline{u}}$
we conclude that
\begin{eqnarray} \label{EF__Equ_Pol1_Rad} \hspace*{3cm}
- \triangle{\;\!}\left({\underline{x}}^{T} (\triangle{\;\!}
{\underline{u}})\right)\,=\,- {\triangle}^{2}{\;\!}({\underline{x}}^{T}{\underline{u}})
&\,=\,& \lambda{\;\!}\triangle{\;\!}({\underline{x}}^{T}{\underline{u}}
)\quad .
\end{eqnarray} 
Using equation (\ref{BewToPo}), namely, ${\bf B}{\chi}\,=\,
{{\underline{x}}}^{T}{\underline{u}}$ and by taking into account
(\ref{Lapl_Komm-Belt}) we find
\begin{eqnarray} \label{EF__Equ_Pol1_Rad1} \hspace*{4cm}
- {\bf B}({\triangle}^{2}{\;\!}{\chi})
&\,=\,& {\bf B}(\lambda{\;\!}{\triangle}{\;\!}{\chi})\quad .
\end{eqnarray} 
A simple consequence of (\ref{ToPo3}) is that ${\chi}$ satisfies
(\ref{ToPo3}) as well as
${\triangle}^{2}{\;\!}{\chi}$ and ${\triangle}{\;\!}{\chi}$.
So do ${\triangle}^{2}{\;\!}{\chi}$ 
and ${\triangle}{\;\!}{\chi}$
in $D({\bf B}_{\bb R}^{-1})$, which justifies the
application of ${\bf B}_{\bb R}^{-1}$ on  (\ref{EF__Equ_Pol1_Rad1})
(locally in $r> 0\:,r\in  [{\sigma},1]$) and thus,
\begin{eqnarray} \label{EF__Equ_Pol1_Rad2} 
\hspace*{4.5cm}
- {\triangle}^{2}{\;\!}{\chi}
&\,=\,& \lambda{\;\!}{\triangle}{\;\!}{\chi}\,.\quad 
\end{eqnarray} 
The validity of (\ref{EF__Equ_Pol1_Rad2}) for $r=0$  is ensured  again by 
density arguments.\\[.1cm]
Finally we attach boundary conditions for ${\chi}$.
The prescribed homogeneous Dirichlet boundary conditions for a poloidal 
eigenfunction ${\underline{u}}\,=\,{\underline{\mathfrak{p}}}$ on
$\partial{\Omega}_{\sigma}$ 
supply the boundary conditions for the correspondent poloidal scalar
${\chi}$. 
By the evaluation of the component in the
${\underline{\mathfrak{e}}}_{r}$ direction we find 
 ${\chi}\,=\,0 $ on $\partial{\Omega}_{\sigma}$ since
in ${\bf B}_{\bb R}(r^{-1}{\chi})\,\in\,D({\bf B}_{\bb R}^{-1})$ for $r=1$ (and $r={\sigma}>0$).
To gather additionally $\frac{\partial{\chi}}{\partial r}\,=\,0 $ on $\partial{\Omega}_{\sigma}$
from ${u}_{\vartheta}\,=\,{u}_{\varphi}\,=\,0$ on $\partial{\Omega}_{\sigma}$
we use there the relation \\[-.4cm]
\begin{eqnarray} \label{EF__Equ_Pol1_Bound} 
{\underline{0}}\,=\,(0,{u}_{\varphi},-{u}_{\vartheta})^{T}\,=\,{\underline{u}}\,\times\,  
{\underline{\mathfrak{e}}}_{r}\,=\,
grad{\;\!}(\frac{\partial(r{\;\!}{\chi})}{\partial r})
\,\times\, {\underline{\mathfrak{e}}}_{r}\quad .
\end{eqnarray}
No we apply $\frac{1}{|\omega_{r}|}\intq_{\omega_{r}}
(\frac{\partial(r{\;\!}{\chi})}{\partial r})({\underline{x}})d
\omega_{r}\,=\,0$ $\forall$ $r\,\in\,[\sigma ,1], \sigma > 0$ which is
an implication of  (\ref{ToPo3})
and use step by step the arguments for  ${\psi}\,=\,0$ on
$\partial{\Omega}_{\sigma}$ above 
and  ${\chi}\,=\,0 $ on $\partial{\Omega}_{\sigma}$. 
In the following problem we combine the eigenvalue problem for the scalar potential $\chi$ of the correspondent poloidal Stokes eigenfunction ${\underline{u}}\,=\,{\underline{\mathfrak{p}}}\,=\,curl{\;\!}curl{\;\!}({\chi{\;\!}\underline{x}})$
(the eigenvalues $\lambda$ coincide again).
\begin{problem} \label{Prob4} We seek solutions  $\chi$  and  $\lambda$ 
(for ${\sigma}:0\leq{\sigma}<1$)  fulfilling:\\[.25cm]
$- \triangle ^{2}{\;\!}{\chi}\,
 = \,\lambda {\;\!}\triangle {\;\!}{\chi}\,\, $ {in} ${\Omega}_{\sigma}$ , with 
$\chi\,=\,\frac{\partial{\chi}}{\partial r}\,=\,0 $ on $\partial{\Omega}_{\sigma}
\quad$
and (cf. (\ref{ToPo3})):\\[.2cm]
\hspace*{1.3cm}
$\frac{1}{|\omega_{r}|}\intq_{\omega_{r}} f \cdot {\chi}({\underline{x}})d \omega_{r}$	 both
$\forall r >0:r\,\in\,[\sigma,1]$ and
$\forall \,f\,=\,f(r)$ as purely radial function.
\end{problem}
\noindent The best way to solve Problem \ref{Prob4} is 
to regard the partial differential equation as a system of two partial
differential equations.
\subsection{The Stokes eigenfunctions in the unit ball \label{sec_sto_ew2}}
We regard the eigenvalue problem of the Stokes Operator ${\bf A}_{o}$
formulated in Section \ref{sec_int} 
as  Problem \ref{Prob1}  for $[{\underline{u}},p] \,\in\,{\underline{\bb C}}^{2}({\Omega}_{o}) \times  {{\bb C}}^{1}({\Omega}_{o})$
and as Problem \ref{Prob2} for 
${\underline{u}}\,\in\,{\underline{\bb S}}^{2}({\Omega}_{o})$, respectively.
Let the  Bessel functions of first kind be denoted by
\begin{eqnarray} 
\label{Bessel}
J_{l+\frac{1}{2}}(t) \,\,, \,J_{-(l+\frac{1}{2})}(t) &{~}&
\, \forall\,l\,\in\,{\bb N}_{o}\,.
\end{eqnarray}
The Bessel function $J_{l+\frac{1}{2}}(t)$ has infinitely many positive zeros 
$\{{\mu}_{l+\frac{1}{2}}^{j}\}_{j=1}^{\infty}$, with $0\,<  {\mu}_{l+\frac{1}{2}}^{1}\,<{\mu}_{l+\frac{1}{2}}^{2},\,< \dots$
and ${\mu}_{l+\frac{1}{2}}^{j}\rightarrow\,\infty $ for
${j}\rightarrow\,\infty $.
\\[3mm]
Applying the normed surface spherical harmonics \eqref{SuSphZ}
$
\{Z_{l}^{-l},\dots,Z_{l}^{-1},Z_{l}^{0},Z_{l}^{1},\dots,Z_{l}^{l}\}
_{l=0}^{\infty}
$ (cf. Appendix {{{{Theorem B}}}) 
we 
define the  following toroidal fields for  $l\,\in\,{\bb N}$ and $k=0,\dots,l$  (resp. $k=1,\dots,l$)
\begin{eqnarray} 
\label{EigTo}
{\underline{\mathfrak{t}}_{c,}}^{(j),l,k}
_{\mathfrak{s}} :=\,(curl( {\frac{J_{l+\frac{1}{2}}( {\mu}_{l+\frac{1}{2}}^{j} r) }{\sqrt{r }}} Z_{l}^{k}   {\;\!}
{\underline{x}}))_{\mathfrak{s}}  & ,&
{\underline{\mathfrak{t}}_{s,}}^{(j),l,k}
_{\mathfrak{s}} :=\,(curl({\frac{J_{l+\frac{1}{2}}( {\mu}_{l+\frac{1}{2}}^{j} r) }{\sqrt{r }}}Z_{l}^{-k}   {\;\!}
{\underline{x}}))_{\mathfrak{s}}
\, \forall\,j,\, l\,\in\,{\bb N}\,.
\end{eqnarray}\\[-.3cm]
For the definition of the poloidal fields we will use the following scalar functions ${\chi}$: $\forall\,j,\, l\,\in\,{\bb N}$ 
\begin{eqnarray} \hspace*{-1.2cm}
{\underline{\mathfrak{p}}_{c,}}^{(j),l,k}_{\mathfrak{s}} &:=& \,(curl(curl({\chi}_{c}^{(j),l,k}{\underline{x}})))_{\mathfrak{s}}\,\,,\,
{\chi}_{c}^{(j),l,k}
 :=\,({\frac{J_{l+\frac{1}{2}}( {\mu}_{l+\frac{3}{2}}^{j} r) }{\sqrt{r }}} -{\frac
{{\mu}_{l+\frac{3}{2}}^{j}J_{l-\frac{1}{2}}( {\mu}_{l+\frac{3}{2}}^{j}) r^{l}}{2l+1}})
Z_{l}^{k}   \,,
  \label{EigPoChi} 
  \\
{\underline{\mathfrak{p}}_{s,}}^{(j),l,k}_{\mathfrak{s}} &:= &\,(curl(curl({\chi}_{s}^{(j),l,k}{\underline{x}})))_{\mathfrak{s}}\,\,,\,
{\chi}_{s}^{(j),l,k}
 :=\,({\frac{J_{l+\frac{1}{2}}( {\mu}_{l+\frac{3}{2}}^{j} r) }{\sqrt{r }}} -{\frac
{{\mu}_{l+\frac{3}{2}}^{j}J_{l-\frac{1}{2}}( {\mu}_{l+\frac{3}{2}}^{j}) r^{l}}{2l+1}})
Z_{l}^{-k} \nonumber
\end{eqnarray}
both for $k=1,\dots,l$.
We collect the toroidal fields and the poloidal fields in the following system of solenodal functions
\begin{eqnarray} 
\label{EigFu}
\{\{{\underline{\mathfrak{p}}_{s,}}^{(j),l,l}_{\mathfrak{s}},\dots, {\underline{\mathfrak{p}}_{c,}}^{(j),l,0}_{\mathfrak{s}}
,\dots,{\underline{\mathfrak{p}}_{c,}}^{(j),l,l}_{\mathfrak{s}}\}\cup \{{\underline{\mathfrak{t}}_{s,}}^{(j),l,l}_{\mathfrak{s}},\dots, {\underline{\mathfrak{t}}_{c,}}^{(j),l,0}_{\mathfrak{s}}
,\dots,{\underline{\mathfrak{t}}_{c,}}^{(j),l,l}_{\mathfrak{s}}\}\}_{j,l\,\in\,{\bb N}}\,.
\end{eqnarray}
Now we are able to formulate our main result. We use the Bessel functions (\ref{Bessel}) and the definitions 
(\ref{EigTo}) and (\ref{EigPoChi}).
\begin{theorem} \label{EFBall_Thm} The solenoidal fields (\ref{EigFu}) form a complete system of eigenfunctions for the Stokes operator 
${\bf A}_{o}$ in the unit ball.
The eigenvalues are  $\lambda_{(j),l,k}\,=\,({\mu}_{l+\frac{1}{2}}^{j})^{2} $  with the multiplicity $(2l+1)$ in the class of 
toroidal fields and the values 
$\lambda_{(j),l,k}\,=\,({\mu}_{l+\frac{3}{2}}^{j})^{2} $ with the
multiplicity  $(2l+1)$ in the class of the poloidal fields.
\end{theorem}
\begin{proof} The fields of the system  (\ref{EigFu}) are solenoidal
  by their definition.  Additionally we see their orthogo\-na\-li\-ty
  through the application of Theorem \ref{Dec_Thm}.  
From the definition of the  ${\underline{\mathfrak{t}}}_{\mathfrak{s}}$
(\ref{Toexpl}) and the  ${\underline{\mathfrak{p}}}_{\mathfrak{s}}$
(\ref{Poexpl}), it is apparent that from the Stokes eigenvalue problem
boundary value problems one obtains ordinary  
differential equations for the functions  ${\hat{\psi}}(r)$ and  ${\hat{\chi}}(r)$ as the Fourier coefficients of  ${\psi}$ and ${\chi}$ developed in the 
surface spherical harmonics. Here the tools of {{Remark A}} in the
appendix are crucial.
The proof is finished by the solution and a case discussion for solutions of the boundary value problems of ordinary differential equations
for both  ${\hat{\psi}}(r)$ and  ${\hat{\chi}}(r)$.\hfill 
\end{proof}

\subsection{The Stokes eigenfunctions in open 3D-Annuli \label{sec_sto_ewAnn}}
In this subsection we study the eigenvalue problem of the Stokes Operator ${\bf A}_{\sigma}$ formulated in section \ref{sec_int}
as the Problem \ref{Prob1}  for $[{\underline{u}},p] \,\in\,{\underline{\bb C}}^{2}({\Omega}_{\sigma}) \times  {{\bb C}}^{1}({\Omega}_{\sigma})$
and as the Problem \ref{Prob2} for 
${\underline{u}}\,\in\,{\underline{\bb S}}^{2}({\Omega}_{\sigma})$,
respectively. We use the methods explained in \cite{Zand} to 
describe the Stokes eigenfunctions for  ${\bf A}_{\sigma}$ and to show the equations for the determination of 
the eigenvalues. Again, 
the Bessel functions of first kind are written again as
\begin{eqnarray} 
\label{Bessel}
J_{l+\frac{1}{2}}(t) \,\,, \,J_{-(l+\frac{1}{2})}(t) &{~}&
\, \forall\,l\,\in\,{\bb N}_{o}\,.
\end{eqnarray}
Applying the normed surface spherical harmonics (cf. Appendix {{{{Theorem B}}}) 
we 
define the  following toroidal fields for  $l\,\in\,{\bb N}$
\begin{eqnarray} 
\label{EigToAnn}
\hspace*{-1.0cm}
{\underline{\mathfrak{t}}_{\,}}^{l,k}
_{\mathfrak{s}} :=\,(curl{\;\!}({\psi}_{\,}^{l,k}{\;\!}{\underline{x}} )
)_{\mathfrak{s}}  & \mbox{, with} &
{\psi}_{\,}^{l,k}\,:=\,\left(a\cdot {\frac{J_{l+\frac{1}{2}}( {\mu} r) }{\sqrt{r }}} \,+\,b \cdot {\frac{J_{-l-\frac{1}{2}}
( {\mu} r) }{\sqrt{r }}}\right) Z_{l}^{k} \,\,,\, 
\end{eqnarray}
for $ k=-l,\dots,0,\dots,l\,$.
For all $l\,\in\,{\bb N}$ and  $ k=-l,\dots,0,\dots,l $ the poloidal
fields are given through  
\begin{eqnarray} \label{EigPoAnn} 
{\underline{\mathfrak{p}}_{\,}}^{l,k}_{\mathfrak{s}} := \,(curl{\;\!}(curl{\;\!}({\chi}_{\,}^{l,k}{\underline{x}})))_{\mathfrak{s}}\,,
{\chi}_{\,}^{l,k}
 :=\,\big(a\cdot {\frac{J_{l+\frac{1}{2}}( {\mu} r) }{\sqrt{r }}} \,+\,b \cdot {\frac{J_{-l-\frac{1}{2}}
( {\mu} r) }{\sqrt{r }}} \,+\,c \cdot r^{l}\,+\,d \cdot r^{-l-1}
\big) 
Z_{l}^{k} \,.
\end{eqnarray}
For fixed  ${\sigma}:0\leq{\sigma}<1$, the eigenvalues $\lambda\,=\,
{\mu}^2$ are determined by finding ${\mu}$ as roots of the
following trans\-cendental 
equations. These have infinitely many positive zeros like the Bessel
functions (\ref{Bessel}). 
We start with the 
toroidal fields. Considering the vanishing boundary conditions on
$\partial{\Omega}_{\sigma}$ and their form (see equation
\eqref{curlToPo}) we obtain them for all  
$l\,\in\,{\bb N}$ and write them as characteristic equation with determinants
\begin{equation}\label{detATo} 
 \text{det}{\;\!}
 \left[\left(
 \begin{array}{cc} 
 {J_{l+\frac{1}{2}}( {\mu} {\sigma}) }& {J_{-l-\frac{1}{2}}( {\mu} {\sigma})}\\
{J_{l+\frac{1}{2}}( {\mu} ) }&{J_{-l-\frac{1}{2}}( {\mu} )}
 \end{array}
 \right)\right]\, =\,0 \quad ,\,\,l\,\in\,{\bb N} \,.
\end{equation}
For fixed $l\,\in\,{\bb N}$ the (positive) roots of \eqref{detATo} $\{{\mu}^{j}_{l}\}_{j=1}^{\infty}$ are distributed like the zeros of Bessel functions, with $0\,<  {\mu}_{l}^{1}\,<{\mu}_{l}^{2},\,< \dots$
and ${\mu}_{l}^{j}\rightarrow\,\infty $ for ${j}\rightarrow\,\infty $. The rank  of the matrix in \eqref{detATo}
is one for the values $\mu\,=\,{\mu}_{l}^{\tilde{j}}$ at
${\mu}_{l}^{\tilde{j}}\,\in\,\{{\mu}^{j}_{l}\}_{j=1}^{\infty}$.
The constant $b$ depends on $a$ and $a$ can be normalized.
\\[3mm]
The transcendental equations for the
poloidal fields ${\underline{\mathfrak{p}}}$ are also imposed by the
vanishing boundary conditions on  $\partial{\Omega}_{\sigma}$. By
investigating the structure of the poloidal fields
${\underline{\mathfrak{p}}}$ 
and the calculation rules for Bessel functions for the values ${\mu}$ where
$l\,\in\,{\bb N}$ (c.f. \cite{Zand}) we succeed.
 We again write a characteristic equation with determinants
\begin{equation}\label{detAPo} 
 \text{det}{\;\!}
 \left[\left(
 \begin{array}{cccc} 
 {J_{l-\frac{1}{2}}( {\mu} {\sigma}) }& {J_{-l+\frac{1}{2}}( {\mu} {\sigma})} & {\sigma}^{l-\frac{1}{2}} & 0\\
{J_{l-\frac{1}{2}}( {\mu} ) }&{J_{-l+\frac{1}{2}}( {\mu} )}& 1 & 0\\
 {J_{l+\frac{3}{2}}( {\mu} {\sigma}) }& {J_{-l-\frac{3}{2}}( {\mu} {\sigma})} & 0 & {\sigma}^{-l-\frac{3}{2}} \\
{J_{l+\frac{3}{2}}( {\mu} ) }& {J_{-l-\frac{3}{2}}( {\mu} )} & 0 & 1
 \end{array}
 \right)\right]\, =\,0 \quad ,\,\,l\,\in\,{\bb N} \,.
\end{equation}
We can also use that for fixed $l\,\in\,{\bb N}$ the (positive) roots
of \eqref{detAPo} $\{{\mu}^{j}_{l}\}_{j=1}^{\infty}$ are distributed
like the zeros of Bessel functions, with $0\,<
{\mu}_{l}^{1}\,<{\mu}_{l}^{2},\,< \dots$ 
with ${\mu}_{l}^{j}\rightarrow\,\infty $ for ${j}\rightarrow\,\infty
$. The rank  of the matrix in \eqref{detAPo} is $3$ for the values
$\mu\,=\,{\mu}_{l}^{\tilde{j}}$ at 
${\mu}_{l}^{\tilde{j}}\,\in\,\{{\mu}^{j}_{l}\}_{j=1}^{\infty}$. From
this it follows that the scalar	 
coefficient $a,\,b$ and $c$ can be written as functions of the parameter $d$. Finally one uses the parameter $d$
to normalize the eigenfunction field 
${\underline{\mathfrak{p}}_{\,}}^{l,k}_{\mathfrak{s}}$ \eqref{EigPoAnn}.
Here the Stokes eigenvalue are given by
$\lambda\,:=\,\left({\mu}_{l}^{\tilde{j}}\right)$ for  
$k=-l,\dots,0,\dots,l\,$.\\[.2cm]
Note that the Stokes eigenvalues and eigenvectors have to be calculated 
whenever the value of ${\sigma}$ is changed. But from another point
of view there is a direct connection 
between the Stokes eigenvalues and the Stokes eigenvalues in ${\Omega}_{R_i,R_a}\,:=\{{\underline{x}} \in {\bb R}^{3}:\,R_i <||{\underline{x}}||_{\cal{E}}<R_a\,<\,\infty\}$
and ${\Omega}_{\sigma}\,:=\{{\underline{x}} \in {\bb R}^{3}:\,{\sigma}<||{\underline{x}}||_{\cal{E}}<1\}$ for
${\sigma}\,:=\frac{R_i}{R_a}$ (cf. e.g. \cite{RRTZAMM}). In this way 
our tools and investigations do
not only proof the completeness of the systems of Stokes eigenfunctions for  ${\bf A}_{\sigma}$. One can use these 
functions to construct Galerkin schemes for the numerical investigations of
turbulent flows in arbitrary balls and spherical annuli.

\section*{Appendix 1: The Laplace-Beltrami operator, surface spherical
  harmonics and
Bessel functions\label{app_Belt}}
Let $\{(r,\vartheta,\varphi)^{T},\,\mbox{with}\,0\leq r\,,\,0\leq\vartheta\leq\pi\,,\,0\leq\varphi < 2\pi\}$ be the usual spherical polar coordinates in ${\bb R}^{3}$.
We introduce the surface spherical harmonics in the way, that simplifies the 
application of the recursion formulas in the construction (cf. Magnus, Oberhettinger) of  
vector fields in surface spherical harmonics.\\[.1cm]
{{\textbf{Notation A:}}}
The associated Legendre functions for $l=0,1,2,\dots$, and
 $ k \in \{0,1,\dots,l\}$ are 
\begin{eqnarray}\label{LegFunc}
P_{l}^{k}(t)\,:=\,\frac{(-1)^{l+k}}{2^{l}l!}
(1-t^{2})^\frac{k}{2}\frac{d^{l+k}}{dt^{l+k}}(1-t^{2})^{l},\,t\,\in\,[-1,1]\,.
\end{eqnarray}
{{\textbf{Notation B:}}} The real valued surface spherical harmonics of 
degree $l\,\in\,{\bb N}_{o}$ on the unit sphere $\omega$ are called 
$C_{l}^{k}$ and $S_{l}^{k}$. They are defined by:
\begin{eqnarray}\label{SuSphHa1}
C_{l}^{k}(\varphi,\vartheta)&\,:=\,&\cos (k\varphi)P_{l}^{k}(\cos\vartheta),
\,l=0,1,2,\dots\,,
k \,\in \,\{0,1,\dots,l\}\,\nonumber\\[-.2cm]
{~} & {~} &\hspace*{10cm} \mbox{and}\\[-.2cm]
S_{l}^{k}(\varphi,\vartheta)&\,:=\,&\sin (k\varphi)P_{l}^{k}(\cos\vartheta),
\,l=0,1,2,\dots\,,
k \,\in \,\{1,2,\dots,l\}\,\,\nonumber\,.
\end{eqnarray}
We explain the {\it{(in ${{\bb L}_{2}(\omega)}$ sense)}}
normed system of surface spherical harmo\-nics by:\\[-.5cm]
\begin{eqnarray}\label{SuSphHa}
\hspace*{-1cm}{{Z}}_{l}^{k}(.,.):=
\sqrt{\frac{\displaystyle{(l-k)!(2l+1)}}{\displaystyle{2^{\delta_{k,0}}(l+k)!2\pi}}}C_{l}^{k}(.,.)
&\mbox{and}&{{Z}}_{l}^{-k}(.,.):=
\sqrt{\frac{\displaystyle{(l-k)!(2l+1)}}{\displaystyle{(l+k)!2\pi}}}
S_{l}^{k}(.,.)\,.
\end{eqnarray}
{{\textbf{Notation C:}}}
The Bessel functions of the first kind of order $l+\frac{1}{2}$ are 
 $J_{l+\frac{1}{2}}(r)$ for $l\in{\bf Z}$.\\[.2cm]
The following theorems collect results for 
the Laplace-Beltrami o\-pe\-ra\-tor (introduced in
Definition \ref{D4}) from the book of Triebel (\cite{Triebel}).
\\[3mm]
{{\textbf{Theorem A:}}} {\em The Laplace-Beltrami o\-pe\-ra\-tor ${\bf B^{o}}$ is 
essentially self-adjoint in ${\bb L}_{2}(\omega)$. The
associated bilinear form 
$\beta(Y,Z)\,:=\,({\bf B^{o}}Y,Z)_{{\bb L}_{2}(\omega)}\,\,,\,\forall \, Y,Z\,\in\,
 D({\bf B^{o}})$ is nonnegative, i.e. 
$\beta(Y,Y)\,\geq\,0\,\forall \, Y\,\in\,
 D({\bf B^{o}})$ and  $\beta$ is 
 ${\bb W}_{2}^{1}(\omega)$-coercive.
The Friedrichs' extension ${\bf B}$ of ${\bf B^{o}}$, with ${\bf B}\,: \,{\bb W}_{2}^{1}(\omega) \,\rightarrow {\bb L}_{2}(\omega)$, is an operator
with pure point spectrum.}\\[.22cm]
{{\textbf{Theorem B:}}} {\em The eigenvalues of ${\bf B}$ are the $\nu_{(l)}:=l(l+1)$ and
the multiplicity
of every $\nu_{(l)}$ is $N(\nu_{(l)}):=2l+1$
$\forall \, l$.
The system ${\bf B}$:
$\{\tilde{Y}_{l}^{l},\dots,\tilde{Y}_{l}^{1},Y_{l}^{0},Y_{l}^{1},\dots,Y_{l}^{l}\}
_{l=0}^{\infty}\,:=\,
\{{{Z}}_{l}^{-l},\dots,{{Z}}_{l}^{-1},
{{Z}}_{l}^{0},{{Z}}_{l}^{1},\dots,
{{Z}}_{l}^{l}\}
_{l=0}^{\infty}\,$ of
eigenfunctions of
is a complete orthonormal system in ${\bb L}_{2}(\omega)$.}\\[.3cm] 
{{\textbf{Theorem C:}}} 
{\em The eigenvalues of ${\bf B}$ can be ordered by their (absolute) values, taking into account their multiplicities. If the $\{\nu_{j}\}_{j=0}^{\infty}$
are the ordered eigenvalues and the $\{Y_{j}\}_{j=0}^{\infty}$ the correspondent 
orthonormal eigenfunctions, then the system  $\{Y_{j}\}_{j=0}^{\infty}$ is 
complete in ${\bb L}_{2}(\omega)$ and 
$D({\bf B})\,=\,\{Y\,\in\,{\bb L}_{2}(\omega):\,
\sum_{j=0}^{\infty}\,(1+\nu_{j}^{2})|(Y,Y_{j})_{{\bb L}_{2}(\omega)}|^{2}
<\infty\}$ with 
${\bf B}Y\,=\,\sum_{j=1}^{\infty}\,
\nu_{j}(Y,Y_{j})_{{\bb L}_{2}(\omega)}Y_{j}$, since 
$\nu_{0}\,=\,0$.}\\[.1cm] 
{{\textbf{Notation D:}}} 
Following \cite[Sec.2.4.1]{Linda} and \cite[Sec. 6.4.5]{Triebel}
we use, that the sets of all harmonic functions on ${\Omega}_{o}$
and on ${\Omega}_{\sigma}\,,\,\sigma > 0,$ respectively, are spanned by the systems
\begin{eqnarray}\label{HarmFunc} 
H_{o}&:=&\quad\{r^{l}\cdot C_{l}^{k}(\vartheta,\varphi)\}^{l}_{k=0,\,l \in {\bb N}_{0}}
\cup \{r^{l}\cdot S_{l}^{k}(\vartheta,\varphi)\}^{l}_{k=1,\,l \in {\bb N}}\hspace*{2cm} \mbox{and} \nonumber\\[-.3cm]
\,&\, \\[-.2cm]
H_{\sigma}&:=&\{\{r^{l}\cdot C_{l}^{k}\}^{l}_{k=0}
\cup 
\{r^{-(l+1)}\cdot C_{l}^{k}\}^{l}_{k=0}\}_{l \in {\bb N}_{0}}
\cup \{\{r^{l}\cdot S_{l}^{k}\}^{l}_{k=1}
\cup \{r^{-(l+1)}\cdot S_{l}^{k}\}^{l}_{k=1}\}_{l \in {\bb N}}\,\,\,,\nonumber
\end{eqnarray}\\[-.6cm]
where the $C_{l}^{k}(\vartheta,\varphi)$ and $S_{l}^{k}(\vartheta,\varphi)$ are the
surface spherical harmonics on the unit sphere
$\omega$ 	in accordance \\
with {{Notation B}}.

\section*{Appendix 2: Auxiliary material and Nabla-Calculus\label{aops_nab}}
We suppose that ${\eta}({\underline{x}})$ and ${\underline{u}}({\underline{x}})$ 
are sufficiently smooth and ${\underline{x}}$ and ${\underline{u}}$ are in ${\bb R}^{3}$.
The following rules are independent of the special choice of the
orthonormal right-handed 
coordinate system in ${\bb R}^{3}$, where ${\underline{\nabla}}$ denotes
the Nabla operator  (cf. \eqref{nab..def:} and Remark \ref{R2}).\\[.2cm]
{{\textbf{Remark A:}}} (Product formulas in the  Nabla-calculus)\\[-.7cm]
\begin{eqnarray}\label{rotgrad}
(i) \, \quad curl\,grad \,{\eta}\, =   \,{\underline{\nabla}}\,
\,\times\,({\underline{\nabla}}\,{\eta})\, =  \underline{0} \quad & {~} &
(ii) \,\quad div\,curl \,{\underline{u}}\, =   \,
{\underline{\nabla}}^{T}\cdot ({\underline{\nabla}}\,\,\times\,
{\underline{u}})\,=  {0}\,\nonumber\\[-.25cm]
{~}& {~} &{~}\\[-.25cm]
(iii) \,\quad curl\,curl\,{\underline{u}}\, = 
{\underline{\nabla}}\,\,\times\, ({\underline{\nabla}}\,\,\times\,
{\underline{u}})\, & = &
grad \,div \,{\underline{u}}\,\,-
\,\triangle\,{\underline{u}} \,  \,.\nonumber
\end{eqnarray} 
{{\textbf{Remark B:}}} (Formulas for products of scalar and
vectorial functions)\\[-.7cm]
\begin{eqnarray}\label{nab_prod}
(i)\,\quad curl \,({\eta}{\underline{u}})\, =  
\,{\underline{\nabla}}\,\,\times\, ({\eta}{\underline{u}})& = & {\eta} \,curl \,{\underline{u}}
\,+\,\,\,(grad \,{\eta}) \,\times\, {\underline{u}}
\quad 
\nonumber\\[-.25cm]
{~}& {~} &{~}\\[-.25cm]
(ii) \quad div \,({\eta}{\underline{u}})\, =  \,{\underline{\nabla}}^{T}\cdot ({\eta}
{\underline{u}})\,& = & 
\,{\eta} \,div \,{\underline{u}}
\,\,+\,\,{\underline{u}}^{T}\cdot\,grad \,{\eta}\, .
\nonumber
\end{eqnarray}
It is advantageous for our considerations to point out relations of
the scalar Laplacian and 
the Laplace-Beltrami operator both regarded in the usual spherical
polar coordinates in ${\bb R}^{3}$ as follows.\\[.1cm]
{{\textbf{Remark C:}}} (Laplacian and the Laplace-Beltrami operator)\,
The scalar Laplacian $\triangle$ applied on the function ${\eta}$ has the form:\\[-.3cm]
\begin{eqnarray}\label{Lapl_Stand}
\triangle{\;\!} {\eta}\,:=\,
{\underline{\nabla}}^{T}\cdot{\underline{\nabla}} {\;\!}{\eta}& = &\frac{\displaystyle{1}}{\displaystyle{r^{2}}}
\frac{\partial}{\partial r}
(r^{2}\frac{\partial {\eta}}{\partial r})\,+\,
\frac{\displaystyle{1}}
{\displaystyle{r^{2}\sin{\theta}}}
\frac{\partial }{\partial \theta}(\sin{\theta}\frac{\partial {\eta}}{\partial \theta})
\,+\,
\frac{\displaystyle{1}}{\displaystyle{r^{2}\sin ^{2}{\theta}}}
\frac{\displaystyle{\partial ^{2} {\eta}}}
{\displaystyle{\partial \phi^{2}}}
\nonumber\\[-.25cm]
{~}& {~} &{~}\\[-.25cm]
{~}& = &\frac{\displaystyle{1}}{\displaystyle{r^2}}
\frac{\partial}{\partial r}\big(r^2\frac{\partial \eta}{\partial r}\big)
\,-\,\frac{\displaystyle{1}}{\displaystyle{r^{2}}}{\bf B}{\;\!} {\eta}\quad\,,\nonumber
\end{eqnarray}
where $\bf B$ denotes the Laplace-Beltrami operator (cf. Definition \ref{D4}  ). 
One shows by a simple calculation that the Laplace-Beltrami operator $\bf B$
and the scalar Laplacian commute for all functions 
$\,{\eta}\,\in\,{\bb C}^{4}({\overline{\Omega}}_{\sigma})$ at least in all inner points ${\underline{x}}$ of ${\overline{\Omega}}_{\sigma}$
(in ${\cal C}^{4}({\Omega}_{\sigma})$) for ${\sigma}:0\leq{\sigma}<1$
\begin{eqnarray}\label{Lapl_Komm-Belt}
\hspace*{4cm}
\triangle{\;\!}{\bf B}{\;\!} {\eta} \,=\,{\bf B}{\;\!} \triangle{\;\!}{\;\!} {\eta}\,& { ~}&\,\forall\,\,{\eta}\,\in\,{\bb C}^{4}({\overline{\Omega}}_{\sigma})\,.
\end{eqnarray}
Next we collect commutative laws
for differential operators, which we need for the 
proofs of the Theorems \ref{Dec_Thm} and \ref{EFBall_Thm}.\\[.2cm]
{{\textbf{Remark D:}}}
Firstly, we regard  an arbitrary vector field of the form:
${\underline{\mathfrak{y}}}\,=\,{\eta}({\underline{x}}){\underline{x}}$.
Without any loss of generality, we use canonical 
(Cartesian) coordinates here. But the derived relations are independent of the 
actual chosen orthonormal right-handed coordinate system, since the definitions 
of the differential operators $grad\,,div\,$ and $curl$
are independent of the chosen orthonormal right-handed coordinate system
too. 
Let us start with the field: $curl{\;\!}(curl{\;\!}{\underline{\mathfrak{y}}})$.
Here we find 
by the use of (\ref{rotgrad}), (\ref{nab_prod}) and simple
calculations that
\begin{eqnarray}\label{Pol-Feld-Trick}
curl{\;\!}(curl{\;\!}({\eta}({\underline{x}}){\underline{x}}))&\,=\,&
grad{\;\!}({\eta}\,+\,{\underline{x}}^{T}\cdot grad{\;\!}{\eta})
\,-\,(\triangle{\;\!}{\eta}) {\underline{x}}\quad 
\end{eqnarray}
(written in spherical polar coordinates). This is 
the justification for (\ref{Poexpl}). 
The repeated application of the $curl$ yields together with (\ref{rotgrad}), (\ref{nab_prod}), and
(\ref{Pol-Feld-Trick}) yield
\begin{eqnarray} \label{TP_Ver_Lap1} 
\triangle{\;\!}(curl{\;\!}({\eta}{\underline{x}}))
&\,=\,& curl{\;\!}((\triangle{\;\!}{\eta}){\underline{x}}))\quad \quad \text{and}
\end{eqnarray} 
\begin{eqnarray}\label{TP_Ver_Lap2} 
\triangle{\;\!}(curl{\;\!}curl{\;\!}({\eta}{\underline{x}}))
&\,=\,& curl{\;\!}curl{\;\!}((\triangle{\;\!}{\eta}){\underline{x}})\quad.
 \end{eqnarray}
One can consider the  properties of toroidal fields and poloidal fields 
(({\ref{curlToPo}}) resp. ({\ref{curl2ToPo}})) as simple conclusions of
the relations (\ref{Pol-Feld-Trick}),
(\ref{TP_Ver_Lap1}) and  (\ref{TP_Ver_Lap2}).\\[.2cm]
{{\textbf{Remark E:}}}
We suppose now, that ${\underline{u}}$ is an  arbitrary solenoidal vector field 
(cf. Notation \ref{N9}).
Since $div{\;\!}{\underline{u}}\,=\,0$ , applying
(\ref{rotgrad}) $(iii)$  we find 
\begin{eqnarray} \label{EF_Ver_Lap1} \hspace*{4cm}
\,-\,\triangle{\;\!}
{\underline{u}}({\underline{x}})
&\,=\,& curl{\;\!}( curl{\;\!} {\underline{u}}({\underline{x}})  )\quad ,
\end{eqnarray} \\[-.4cm] 
and furthermore the relation:\\[-.7cm]
\begin{eqnarray} \label{EF_Ver_Lap2} \hspace*{4cm}
\triangle{\;\!}(\underline{x}^{T}
{\underline{u}})\,=\,{\underline{x}}^{T}(\triangle{\;\!}{\underline{u}})
\,+\,2{\;\!}div{\;\!}{\underline{u}}
&\,=\,&{\underline{x}}^{T}(\triangle{\;\!}{\underline{u}})\quad .
\end{eqnarray}  	 
Using these relations we receive for the radial components of $curl{\;\!} (\triangle{\;\!}{\underline{u}})$ and 
$curl{\;\!} curl{\;\!}(\triangle{\;\!}{\underline{u}})$
\begin{eqnarray} \label{EF_Tol_Lap1} 
{\underline{x}}^{T}curl{\;\!}\triangle{\;\!}{\underline{u}}\,=\,
{\underline{x}}^{T}(\triangle{\;\!}curl{\;\!}{\underline{u}})
&\,=\,& \triangle{\;\!}({\underline{x})^{T} (curl{\;\!} 
{\underline{u}}))}
\end{eqnarray} 
and
\begin{eqnarray} \label{EF_Pol_Lap2}
-{\underline{x}}^{T}  (curl{\;\!} curl{\;\!} (\triangle{\;\!}{\underline{u}}))\,=\,
-{\underline{x}}^{T}  \triangle{\;\!}(curl{\;\!} curl{\;\!} {\underline{u}})
&\,=\,&\triangle{\;\!}({\underline{x}}^{T}(\triangle{\;\!}{\underline{u}}))\,=\,\triangle^{2}{\;\!}({\underline{x}}^{T}
{\underline{u}})\,.
\end{eqnarray}


\end{document}